\documentclass[onefignum,onetabnum,onealgnum,oneeqnum,reqno]{siamart250211}
\usepackage[papersize={8.5in,11in},margin=1.0in,footskip=0.3in]{geometry}
\usepackage[final]{microtype}
\usepackage{amssymb,amsmath}
\usepackage{bm}
\usepackage{enumitem}
\usepackage{mathrsfs}
\usepackage{booktabs}
\usepackage{multirow}
\usepackage{tcolorbox}
\usepackage{tikz}
\usepackage{orcidlink}

\setlist[enumerate]{leftmargin=2.6em,label=(\roman*),topsep=0.5em,parsep=0.25em}
\setlist[itemize]{leftmargin=2em,topsep=0.25em,parsep=0.25em}

\definecolor{bluey}{HTML}{008fba}
\hypersetup{urlcolor=bluey,linkcolor=bluey,citecolor=bluey}

\usepackage{silence}
\usepackage[font={sf,small},labelsep=period,labelfont=bf,margin=6mm,skip=1mm,belowskip=-1mm]{caption}

\definecolor{bluesection}{HTML}{1e4b5e}
\makeatletter
\renewcommand\section{\@startsection{section}{1}{.25in}{1.3ex \@plus .5ex \@minus .2ex}{-.5em \@plus -.1em}{\reset@font\normalsize\bfseries\color{bluesection}}}
\renewcommand\subsection{\@startsection{subsection}{2}{.25in}{1.3ex\@plus .5ex \@minus .2ex}{-.5em \@plus -.1em}{\reset@font\normalsize\bfseries\color{bluesection}}}
\renewcommand\subsubsection{\@startsection{subsubsection}{3}{.25in}{1.3ex\@plus .5ex \@minus .2ex}{-.5em \@plus -.1em}{\reset@font\normalsize\bfseries\color{bluesection}}}
\makeatother

\crefname{figure}{Fig.\!}{Figs.\!}
\crefname{section}{\S\!}{\S\!}
\crefname{subsection}{\S\!\!}{\S\!\!}
\newcommand{\R}{\mathbb R}
\newcommand{\vu}{\mathbf u}
\newcommand{\vn}{\mathbf n}
\newcommand{\vf}{\mathbf f}
\newcommand{\vg}{\mathbf g}
\newcommand{\vh}{\mathbf h}
\newcommand{\vtau}{\bm \tau}
\newcommand{\vsigma}{\bm \sigma}
\newcommand{\trans}{{\mathsf{T}}}
\newcommand{\jump}[1]{[\![ #1 ]\!]}
\newcommand{\polydeg}{\wp}
\DeclareMathOperator{\diag}{diag}
\DeclareMathOperator*{\argmin}{argmin}
\DeclareMathOperator*{\spanop}{span}
\newcommand{\xin}{x}
\newcommand{\xout}{\bar{x}}

\usepackage{algorithm}
\usepackage[noend]{algpseudocode}
\makeatletter
\algrenewcommand\ALG@beginalgorithmic{\footnotesize}
\renewcommand{\ALG@name}{\small Algorithm}
\makeatother

\definecolor{mbluey}{HTML}{52829c}
\newtcbox{\inlinebox}[1][]{on line,boxrule=0.5pt,boxsep=0pt,top=0.3pt,left=1pt,bottom=-0.2pt,right=1pt,colback=white,colframe=mbluey,arc=3pt,#1}
\newcommand{\solver}[1]{\inlinebox{\sffamily#1\vphantom{hg}}}

\newcommand{\SAI}[1]{\ensuremath{\mathsf{SAI\scalebox{0.875}{$\mathsf{#1}$}}}}
\newcommand{\xstack}[2]{\raisebox{-0.4ex}{\scalebox{1.4}[1.9]{$\substack{#1\\[-1.1ex]#2}$}}}
\newcommand{\xOO}{\xstack{\rightharpoonup}{\rightharpoondown}}
\newcommand{\xOI}{\xstack{\rightharpoonup}{\leftharpoondown}}
\newcommand{\xIO}{\xstack{\leftharpoonup}{\rightharpoondown}}
\newcommand{\xII}{\xstack{\leftharpoonup}{\leftharpoondown}}
\newcommand{\xbar}{\raisebox{0.22ex}{$|$}}
\newcommand{\xCG}{\ensuremath{\mathsf{CG}}}
\newcommand{\xGMRES}{\raisebox{0.2ex}{$\substack{\mathsf{GM}\\[-0.4ex]\scriptscriptstyle\mathsf{RES}}$}}
\newcommand{\xGSsymbol}{\includegraphics[height=1.6ex]{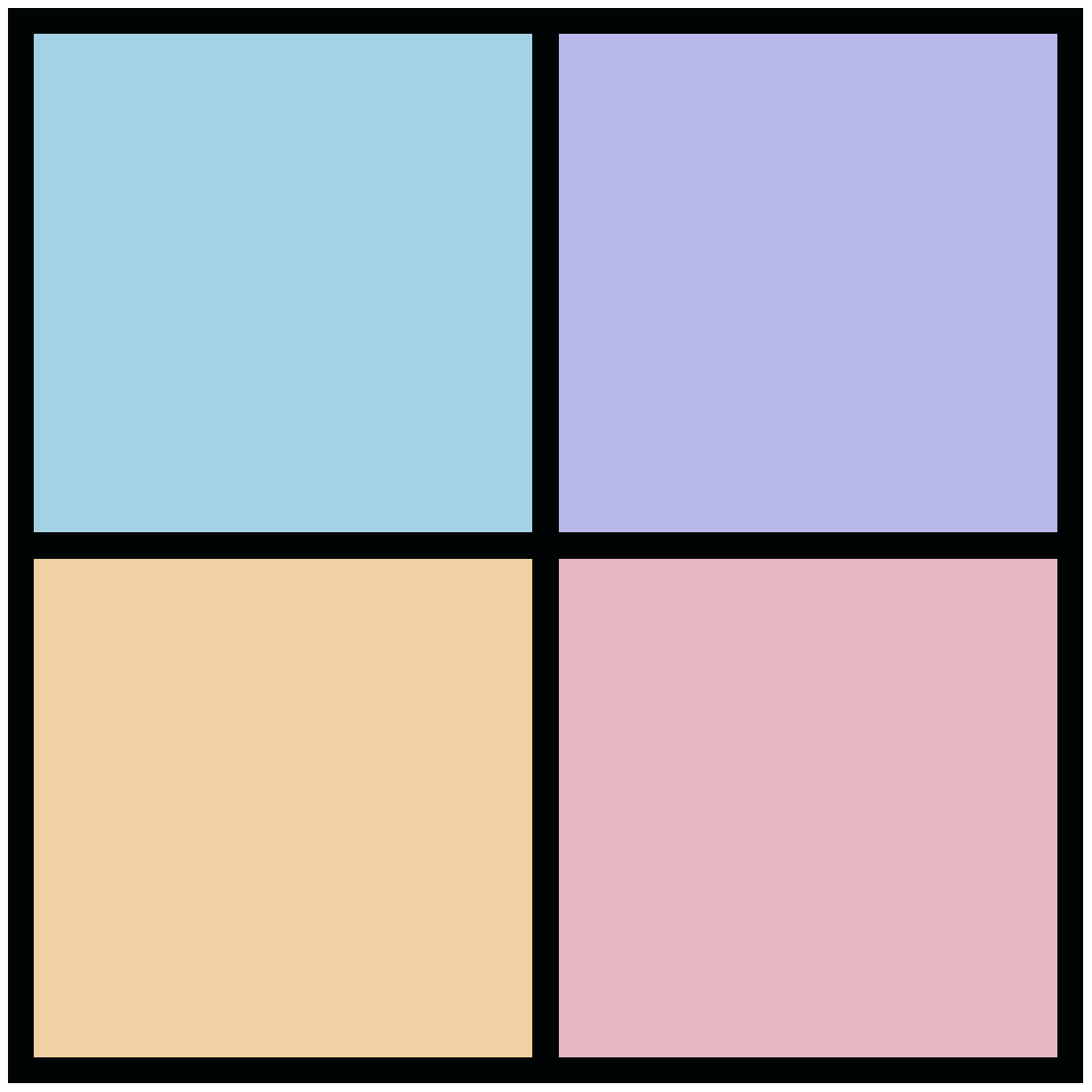}}
\newcommand{\xLGSsymbol}{\includegraphics[height=1.7ex]{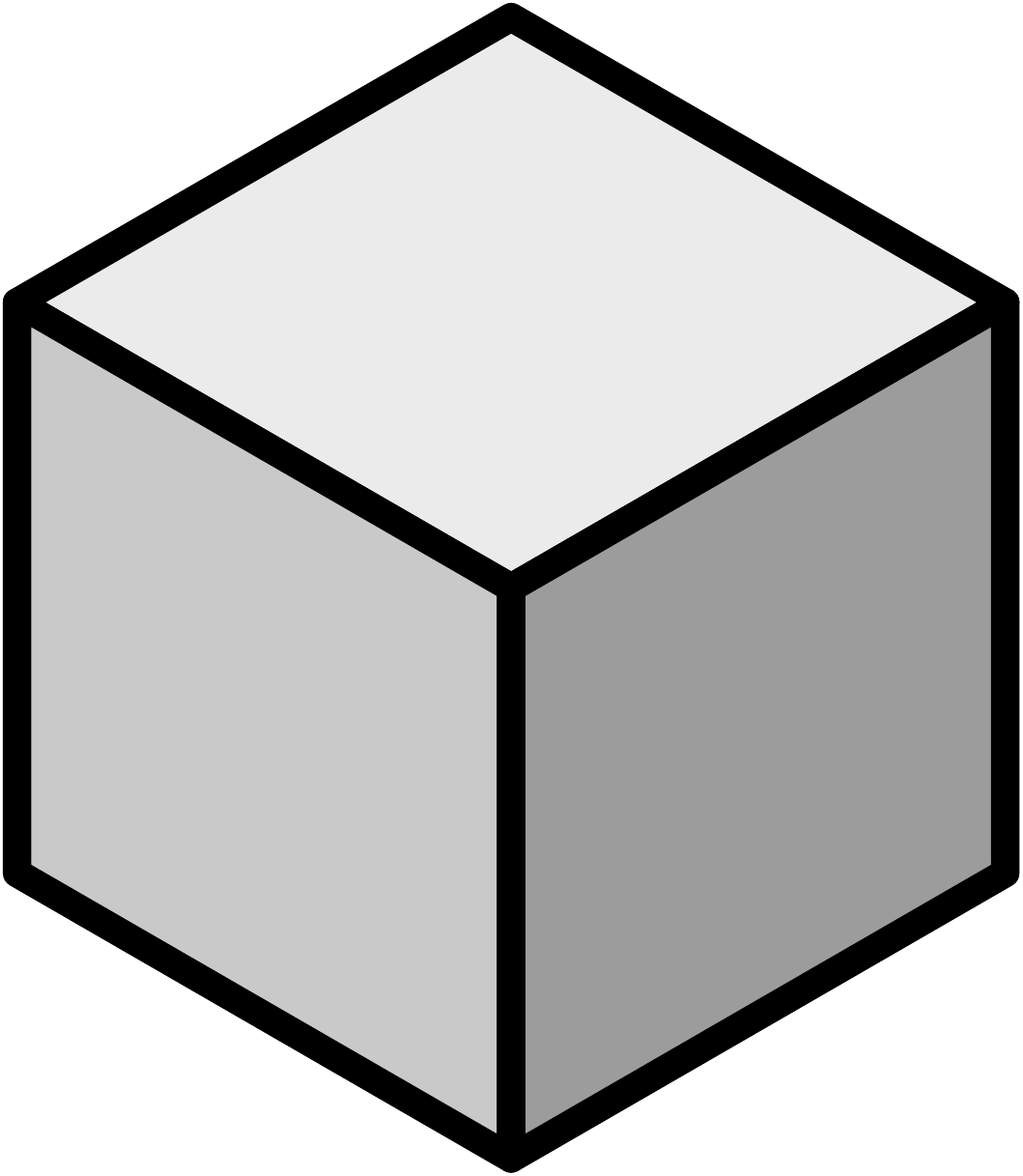}}
\newcommand{\xGS}{\ensuremath{\xGSsymbol\mathsf{GS}}}
\newcommand{\xLGS}{\ensuremath{\xLGSsymbol\mathsf{GS}}}
\newcommand{\xJCG}{\solver{CG\xbar J}}
\newcommand{\xJGL}{\solver{\xGMRES\xbar J}}
\newcommand{\xGSOOGL}{\solver{\xGMRES\xbar\xGS\xOO}}
\newcommand{\xGSOIGL}{\solver{\xGMRES\xbar\xGS\xOI}}
\newcommand{\xGSOICG}{\solver{CG\xbar\xGS\xOI}}
\newcommand{\xLGSOOGL}{\solver{\xGMRES\xbar\xLGS\xOO}}
\newcommand{\xLGSOIGL}{\solver{\xGMRES\xbar\xLGS\xOI}}
\newcommand{\xLGSOICG}{\solver{CG\xbar\xLGS\xOI}}
\newcommand{\xSaiOICG}[1][\ell]{\solver{CG\xbar\SAI{#1}\xOI}}
\newcommand{\xSaiIOCG}[1][\ell]{\solver{CG\xbar\SAI{#1}\xIO}}
\newcommand{\xSaiOOGL}[1][\ell]{\solver{\xGMRES\xbar\SAI{#1}\xOO}}
\newcommand{\xSaiOIGL}[1][\ell]{\solver{\xGMRES\xbar\SAI{#1}\xOI}}
\newcommand{\xSaiIOGL}[1][\ell]{\solver{\xGMRES\xbar\SAI{#1}\xIO}}
\newcommand{\xSaiIIGL}[1][\ell]{\solver{\xGMRES\xbar\SAI{#1}\xII}}

\title{A comparative study of efficient multigrid solvers for high-order local discontinuous Galerkin methods: Poisson, elliptic interface, and multiphase Stokes problems}
\author{Robert I.~Saye \orcidlink{0000-0001-7256-6941}}
\date{\today}

\begin{document}
\setlength{\floatsep}{5mm}
\setlength{\textfloatsep}{5mm}

\makeatletter
\begin{center}
{\bfseries\MakeUppercase\@title\par}%
\vskip2mm%
{\small\textsc{\@author}\par}%
\vskip1mm%
{\sffamily\footnotesize Lawrence Berkeley National Laboratory, Berkeley, California, USA $\mid$ \texttt{rsaye@lbl.gov} \par}%
\vskip5mm%
\end{center}
\makeatother

\begin{abstract}
We design and investigate a variety of multigrid solvers for high-order local discontinuous Galerkin methods applied to elliptic interface and multiphase Stokes problems. Using the template of a standard multigrid V-cycle, we consider a variety of element-wise block smoothers, including Jacobi, multi-coloured Gauss-Seidel, processor-block Gauss-Seidel, and with special interest, smoothers based on sparse approximate inverse (SAI) methods. In particular, we develop SAI methods that: (i) balance the smoothing of velocity and pressure variables in Stokes problems; and (ii) robustly handles high-contrast viscosity coefficients in multiphase problems. Across a broad range of two- and three-dimensional test cases, including Poisson, elliptic interface, steady-state Stokes, and unsteady Stokes problems, we examine a multitude of multigrid smoother and solver combinations. In every case, there is at least one approach that matches the performance of classical geometric multigrid algorithms, e.g., 4 to 8 iterations reduce the residual by 10 orders of magnitude. We also discuss their relative merits with regard to simplicity, robustness, computational cost, and parallelisation.
\end{abstract}

\begin{keywords}
multigrid, local discontinuous Galerkin, elliptic interface, Stokes, multiphase, sparse approximate inverse, Gauss-Seidel
\end{keywords}

\begin{AMS}
65N55 (primary), 65N30, 65N22, 76D07, 35J57
\end{AMS}

\section{Introduction}

In this study, we design and investigate a variety of multigrid methods for solving elliptic interface problems and multiphase Stokes problems discretised through high-order local discontinuous Galerkin (LDG) methods. In particular, using the template of a standard multigrid V-cycle, we consider: (i) an assortment of element-wise block smoothers, including whether to ``reverse'' the smoother after the coarse grid correction; (ii) balancing methods to robustly handle large jumps in ellipticity/viscosity coefficient; and (iii) for Stokes problems, balancing methods to account for the viscosity- and mesh-size dependent scaling of the viscous, pressure gradient, and velocity divergence operators. Many different multigrid smoother and solver combinations are explored here; we numerically examine their performance (i.e., convergence rates) across a wide array of multiphase elliptic interface and Stokes problems, and draw some conclusions about their relative merits.

One particular aim of this study is to explore the potential benefits of sparse approximate inverse (SAI) multigrid smoothers. Given a matrix $A$, SAI methods solve a Frobenius norm minimisation problem (in many cases, $\argmin_B \|AB - I\|_F$) to find the best possible approximate inverse $B$ with a given predefined sparsity structure; for example, a (block) diagonal $B$ would lead to a smoother resembling classical (block) Jacobi, while progressively more dense patterns should result in more accurate inverses and perhaps faster multigrid smoothers. Compared to classical multigrid smoothers, e.g., Jacobi and Gauss-Seidel, SAI-based smoothers could have a number of compelling benefits, including:
\begin{itemize}
    \item \textit{Increased parallelism} \cite{1997doi:10.1137/S1064827594276552,1999HUCKLE1999291,2000doi:10.1137/S106482759833913X,2000doi:10.1137/S0895479899339342,2001doi:10.1137/S1064827500380623,2002BROKER200261}. By virtue of the Frobenius-norm minimisation, each SAI column can be computed independently of every other column, thereby enabling relatively simple and scalable parallelisation approaches. Moreover, for matrices $A$ arising from PDE discretisations on structured meshes, there is often considerable repetition within $A$ that implies many columns of $B$ are identical, up to a one-to-one mapping or reordering of the elements involved; consequently, it is often possible to cache and reuse these solutions across multiple columns of $B$. Put together, it may be the case that one can efficiently compute (in parallel) an SAI whose application as a multigrid smoother yields significantly faster convergence rates than other smoothers.

    Compared to Gauss-Seidel specifically, SAI-based smoothers may also lead to simpler parallel implementations. In particular, Gauss-Seidel usually works best when it is coloured appropriately, e.g., red-black Gauss-Seidel for classical Poisson problems. Multi-coloured Gauss-Seidel methods often lead to highly efficient multigrid solvers, owing to their fast convergence rates and because identically coloured elements can be updated concurrently. On the other hand, they can be complex to implement, especially on large scale unstructured meshes. %
    An alternative approach is processor-block Gauss-Seidel \cite{AdamsBrezinaHuTuminaro} which naturally arises from standard domain decomposition methods: each subdomain is assigned to a processor, with a ghost layer of elemental values synchronised and frozen prior to the parallel execution of Gauss-Seidel confined to each subdomain. In the limit of one element per process, this method decays to (block) Jacobi, which usually needs damping. Consequently, processor-block Gauss-Seidel also needs some amount of damping, the optimal value of which can be subtle to determine especially on unstructured meshes; its convergence can also depend on the runtime processor count. In contrast, SAI-based smoothers usually require no damping, have reproducible behaviour independent of processor/thread count, and only need straightforward (block sparse) matrix-vector products that can be done in parallel.

    \item \textit{Higher arithmetic density}. Computation of an SAI amounts to solving, for each of its (possibly blocked) columns, a relatively small, relatively dense, overdetermined least squares problem. To that end, one can make use of highly optimised dense linear algebra libraries, often implemented via high arithmetic density algorithms (i.e., relatively high floating-point operations per data/memory size), e.g., Level 3 BLAS. For example, although computing a more dense SAI incurs additional floating point operations, some of this may come ``for free'' or with sufficiently minimal overhead compared to the expected faster convergence.

    \item \textit{Increased robustness}. In the context of high-order LDG methods, it is usually necessary to include some kind of penalty stabilisation through user-defined penalty parameters, e.g., to weakly enforce Dirichlet boundary conditions. Oftentimes the performance of Jacobi- and Gauss-Seidel-based multigrid methods can be sensitive to the choice of penalty parameter; although a multigrid-preconditioned Krylov solver can significantly mitigate these effects, it can still be the case that ``badly chosen'' penalty parameters can impact multigrid performance \cite{dgmg,fluxx,flame,FEHN2020109538}. While exploring SAI-based smoothers in this work, we found they can be less sensitive to these effects.

    \item \textit{Adaptivity} \cite{1997doi:10.1137/S1064827594276552,2000doi:10.1137/S106482759833913X,2000doi:10.1137/S0895479899339342,2001doi:10.1137/S1064827500380623,2002BROKER200261,2010doi.org/10.1155/2010/930218}. SAI-based solvers have the potential to automatically and adaptively handle various effects such as anisotropy, highly unstructured meshes, complex geometry, discontinuous PDE coefficients, etc. In particular, a variety of approaches have been developed to automatically refine the sparsity pattern based on some kind of controlled fill-in, e.g., by inserting into the sparsity pattern the next entry or block that best decreases the residual $\|AB - I\|_F$. In the presence of local anisotropy, for example, this approach grows the SAI's effective stencil by prioritising strongly-coupled elements. Fully automated SAI methods could thus be considered as ``black box'' multigrid smoothers.
\end{itemize}
For more background on SAI methods, particularly their use as multigrid smoothers, see, e.g., \cite{2000doi:10.1137/S0895479899339342,2001doi:10.1137/S1064827500380623,2002BROKER200261,2012Sedlacek} and the references therein.

This study builds on our prior work in developing fast multigrid solvers for high-order accurate LDG frameworks \cite{dgmg,fluxx,flame}, as will be elaborated upon in subsequent sections. In particular, the cited works solely considered element-wise block Gauss-Seidel methods, where each block corresponds to the collective degrees of freedom on each mesh element. With a view to the above-noted potential benefits of SAI, we compare these approaches to block SAI-based multigrid solvers. (We also include a comparison to block Jacobi, for interest.)

The paper is structured into two main sections, one for elliptic interface problems, followed by one for multiphase Stokes problems. In each case, we outline the prototype PDE problem and summarise the corresponding LDG framework. In \cref{sec:smoothers}, we define the three essential kinds of multigrid smoother, including a balancing method specifically for SAI in \cref{sec:SAI}, later extended to Stokes problems in \cref{sec:stokesbalance}. Numerical results are presented in \cref{sec:scalarexp} and \cref{sec:stokesexp}, examining a variety of multigrid smoother and solver combinations on an array of increasingly challenging test problems. Key findings are summarised \cref{sec:conclusions} along with a concluding perspective.

\section{Elliptic Interface Problems}
\label{sec:scalar}

The prototype problem we consider here consists of solving for a function $u: \Omega \to \R$ such that
\begin{equation} \label{eq:scalar1} \begin{aligned} -\nabla \cdot \bigl(\mu_i \nabla u \bigr) &= f && \text{in } \Omega_i, \end{aligned} \end{equation}
subject to the interfacial jump conditions
\begin{equation} \label{eq:scalar2} \left. \begin{aligned} \jump{u} &= g_{ij} \\ \jump{(\mu \nabla u) \cdot \vn} &= h_{ij} \end{aligned} \right\} \text{ on } \Gamma_{ij}, \end{equation}
and boundary conditions
\begin{equation} \label{eq:scalar3} \begin{aligned} u &= g_{\partial} && \text{on } \Gamma_D, \\ (\mu \nabla u) \cdot \vn &= h_{\partial} && \text{on } \Gamma_N, \end{aligned} \end{equation}
where $\Omega$ is a domain in $\mathbb R^d$ divided into one or more subdomains (``phases'') $\Omega_i$, $\Gamma_{ij} := \partial \Omega_i \cap \partial \Omega_j$ is the interface between phase $i$ and $j$, and $\Gamma_D$ and $\Gamma_N$ denote the parts of $\partial \Omega$ upon which Dirichlet or Neumann boundary conditions are imposed, respectively. Here, $\jump{\cdot}$ denotes the jump in a quantity across an interface $\Gamma_{ij}$, while $\bm n$ denotes an appropriately-oriented unit normal on $\Gamma_{ij}$ or $\partial \Omega$. Finally, $\mu_i$ is a phase-dependent ellipticity coefficient (``viscosity''), while $f$, $g$, and $h$ provide the data for the elliptic interface problem and are given functions defined on $\Omega$, its boundary, and internal interfaces.

\subsection{LDG schemes}

The LDG framework used in this study is essentially identical to that developed in prior work \cite{dgmg,fluxx}. As such, we provide a summary here and refer the reader to the cited references for details. We consider meshes arising from Cartesian grids, quadtrees or octrees, and adopt a tensor-product polynomial space. Let $\mathcal E = \bigcup_i E_i$ denote the set of mesh elements, let $\polydeg \geq 1$ be an integer, and let $\mathcal Q_\polydeg$ be the space of tensor-product polynomials of one-dimensional degree $\polydeg$; e.g., $\mathcal Q_2$ is the space of biquadratic or triquadratic polynomials in 2D or 3D, respectively. Define the corresponding discontinuous piecewise polynomial space as
\[ V_h = \bigl\{ u : \Omega \to \R \enskip \bigl| \enskip u|_E \in \mathcal Q_{\polydeg} \text{ for every } E \in \mathcal E \bigr\}, \]
so that an arbitrary function of $V_h$ has $(\polydeg+1)^d$ degrees of freedom per element. We similarly define $V_h^d$ and $V_h^{d \times d}$ to be the space of piecewise polynomial vector-valued and matrix-valued fields, respectively.

Solving for a discrete solution $u_h \in V_h$, the high-order LDG method for \cref{eq:scalar1,eq:scalar2,eq:scalar3} is built according to a four-step process: (i) using viscosity-upwinded one-sided numerical fluxes \cite{fluxx}, define a discrete gradient $\bm \tau_h \in V_h^d$ equal to the discretisation of $\nabla u_h$, taking into account Dirichlet source data $g$; (ii) define $\bm \sigma_h \in V_h^d$ as the $L^2$ projection of $\mu \bm \tau_h$ onto $V_h^d$; (iii) using the adjoint of the discrete gradient operator from the first step, compute the negative discrete divergence of $\bm \sigma_h$, taking into account Neumann data $h$, and set the result equal to the $L^2$ projection of $f$ onto $V_h$; finally, (iv) amend the overall system by appropriately including penalty stabilisation. Details of this construction are given in \cite{fluxx}; nevertheless it is worthwhile to highlight some key aspects:
\begin{itemize}
    \item For multiphase problems, we use \textit{viscosity-upwinded} one-sided numerical fluxes to appropriately handle large jumps in ellipticity coefficient $\mu$ \cite{fluxx}. On interfacial element boundaries, situated on the interface between two phases, a viscosity-upwinded flux defines the LDG numerical flux for $\hat{u}$ to source from the more viscous phase (taking into account jump conditions $\jump{u} = g$) and the numerical flux for $\hat{\bm \sigma}$ to source from the less viscous phase (taking into account jump conditions $\jump{(\mu \nabla u) \cdot \vn} = h$). Roughly speaking, in the limit of arbitrarily-large viscosity ratios, this biasing corresponds to the multiphase elliptic interface problem separating into two distinct problems: one for the highly viscous phase (which effectively ``sees'' a Neumann-like boundary condition at the interface), and one for the other phase (which effectively ``sees'' a Dirichlet-like boundary condition). As explored in detail in \cite{fluxx}, biasing the numerical fluxes in this way is key to achieving efficient multigrid solvers as well as optimal high-order accuracy in the LDG discretisation.
    \item Our implementation uses a tensor-product Legendre polynomial modal basis. Nevertheless, we emphasise the algorithms, analysis, and essential conclusions drawn in this work should apply to practically any (reasonable) basis choice, owing to the fact that all of the algorithms operate in element-wise block form. For example, one could also use a Gauss-Lobatto nodal basis.
    \item Unless otherwise stated, we have used the same penalty stabilisation parameters as detailed in prior work \cite{fluxx,flame}.
    \item Assuming sufficiently smooth problem data, optimal order high-order accuracy is attained, i.e., order $\polydeg + 1$ in the maximum norm.
\end{itemize}
The LDG scheme for solving \cref{eq:scalar1,eq:scalar2,eq:scalar3} results in a symmetric positive (semi)definite linear system of the form $A_h u_h = b_h$, where $b_h$ collects the entire influence of the source data $f$, $g$, and $h$ onto the right hand side. Here, $A_h$ is a discrete Laplacian operator and takes the form
\[ A_h = \sum_{i=1}^d G_i^\trans M_\mu G_i + E, \]
where $G = (G_1, \ldots, G_d): V_h \to V_h^d$ is a discrete gradient operator, $M_\mu$ is a $\mu$-weighted block diagonal mass matrix, and $E$ is a block-sparse matrix implementing penalty stabilisation. (We shall also refer to this system as $A_h x_h = b_h$ where $x_h$ relabels $u_h$; the subscript $h$ may be dropped if context permits.)

\subsection{Multigrid schemes}

The multigrid methods considered here apply the standard concepts of a V-cycle operating on a nested mesh hierarchy, combining straightforward interpolation and restriction operators along with a smoothing operator such as block Gauss-Seidel. To summarise some key ingredients:
\begin{itemize}
    \item \textit{Mesh hierarchy.} In this work we use quadtrees and octrees to define the finest mesh; the corresponding tree structure naturally defines a hierarchical process by which to agglomerate fine mesh elements to create a nested mesh hierarchy ${\mathcal E}_h$, ${\mathcal E}_{2h}$, ${\mathcal E}_{4h}$, $\ldots$.
    \item \textit{Interpolation operator.} The interpolation operator $I_{2h}^h : V_{2h} \to V_h$ transfers coarse mesh corrections to a fine mesh and is defined naturally via polynomial injection: $(I_{2h}^h u)|_{E_f} := u|_{E_c}$, where $E_f$ is a fine mesh element and $E_c \supseteq E_f$ is its corresponding coarse mesh element.
    \item \textit{Restriction operator.} The restriction operator $R_h^{2h} : V_h \to V_{2h}$ transfers the residual from a fine mesh to a coarse mesh and is defined via $L^2$ projection, or, equivalently, via the adjoint of the interpolation operator: $R_h^{2h} := M_{2h}^{-1} (I_{2h}^h)^\trans M_h$, where $M_h$ and $M_{2h}$ are the mass matrices on the fine and coarse mesh, respectively.
\end{itemize}
To define the primary operator (i.e., the discrete Laplacian $A_h$) on every level of the hierarchy, there are two possible approaches: (i) one may use a ``pure geometric multigrid'' approach in which every level is explicitly meshed (e.g., listing out the elements, faces, and their connectivity, etc.) upon which the LDG discretisation is explicitly built so as to determine $A_{2h}$, $A_{4h}$, $\ldots$; for best results, the coarse grid problems should apply the same numerical flux and penalty penalisation mechanisms as the fine grid problem; alternatively, (ii) one may use an operator-coarsening method that recursively automates this process. These two methods yield equivalent coarse grid discretisations; however, the latter approach is simpler to implement and automatically guarantees the coarse grid problems are discretised in a way that harmoniously matches the fine grid's discretisation. As such, we have used the latter approach in this work---for details, the reader is referred to \cite{dgmg}, which introduces operator-coarsening methods appropriate to LDG methods, and to \cite{fluxx,flame}, for their extension to variable-viscosity multiphase problems and (un)steady Stokes problems.

\begin{figure}%
\centering%
\begin{minipage}{0.52\textwidth}%
\begin{algorithm}[H]
	\caption{\sffamily\small Multigrid V-cycle $V({\mathcal E}_h, x_h, b_h)$ with $\nu_1$ pre- and $\nu_2$ post-smoothing steps on mesh ${\mathcal E}_h$ of the hierarchy.}
	\begin{algorithmic}[1]
		\If{${\mathcal E}_h$ is the bottom level}
			\State Solve $A_h x_h = b_h$ using the bottom level direct solver.
		\Else
			\State Apply pre-smoother $\nu_1$ times.
			\State Compute restricted residual, $r_{2h} := (I_{2h}^h)^\trans (b_h - {A}_h x_h)$. \label{line:restrict}
			\State Solve coarse grid problem, $x_{2h} := V({\mathcal E}_{2h}, 0, r_{2h})$.
			\State Interpolate and correct, $x_h \leftarrow x_h + I_{2h}^h x_{2h}$.
			\State Apply post-smoother $\nu_2$ times.
		\EndIf
		\State{\textbf{return} $x_h$.}
	\end{algorithmic}%
	\label{algo:vcycle}%
\end{algorithm}%
\end{minipage}%
\end{figure}

We explore here the design of efficient LDG multigrid solvers using the template of a standard V-cycle, as shown in \cref{algo:vcycle}.\footnote{\setlength\lineskiplimit{-\maxdimen}\setlength\baselineskip{1.1\baselineskip}Note that $(I_{2h}^h)^\trans$ appears on line \ref{line:restrict} rather than the operator $R_h^{2h}$; this follows from the Galerkin approach, namely that $A_h x_h = b_h$ implements the PDE problem in the ``mass basis'', i.e., the discretised operator, mapping functions from $V_h$ to $V_h$, is given by $M_h^{-1} A_h$. As a piecewise polynomial function in $V_h$, the residual of the fine mesh problem is given by $M_h^{-1} (b_h - A_h x_h)$. In the multigrid V-cycle, this residual is multiplied by $R_h^{2h}$ to define the right-hand side data for the coarse mesh problem: $M_{2h}^{-1} A_{2h} x_{2h} = R_h^{2h} M_{h}^{-1} (b_h - A_h x_h)$, i.e., $A_{2h} x_{2h} = (I_{2h}^h)^\trans (b_h - A_h x_h) =: r_{2h}$.} The choice of smoother is one of our main focuses. Besides the smoother, we have made the following implementation decisions:
\begin{itemize}
    \item \textit{Bottom solver.} On the bottom level of the multigrid hierarchy, a symmetric eigendecomposition of $A$ is precomputed and used as a direct solver. Whenever the PDE problem at hand has a nontrivial kernel (e.g., constant functions for Poisson problems or constant velocity/pressure fields for Stokes problems, depending on the boundary conditions), a standard pseudoinverse procedure is applied to appropriately handle the kernel eigenmodes.
    \item \textit{Smoother iteration count.} There is a balance between choosing $\nu$ relatively high, thus potentially needing fewer V-cycles, versus choosing $\nu$ relatively low, thus needing more V-cycle iterations. Our prior work in fast LDG multigrid methods has consistently indicated that $\nu_1 \approx 3$ and $\nu_2 \approx 3$ generally yield the best results, based on the metric of fastest computation time for a fixed-factor residual reduction. Consequently, and mainly to help reduce the search space in exploring various smoother choices, throughout this work we have fixed $\nu_1 = \nu_2 = 3$. Further comments on this aspect are given in the concluding remarks.
\end{itemize}
Although multigrid can be used as a standalone iterative method, convergence can be further accelerated by using it as a preconditioner in a Krylov method \cite{dgmg}; two good choices are the Conjugate Gradient (CG) and Generalised Minimal Residual (GMRES) methods \cite{demmel,saad}. To that end, invoking the V-cycle with an initial guess of zero yields a linear operator acting on the given right hand side; this linear operator is denoted as $V$ and is used as the preconditioner for CG and/or GMRES. Note that CG requires a symmetric positive (semi)definite preconditioner, which in turn places some symmetry requirements on the pre- and post-smoothers, discussed further in \cref{sec:scalarexp}. In contrast, GMRES is more flexible and permits a wider array of pre-/post-smoother combinations. Further, with GMRES one could use $V$ as a left-preconditioner (solving $V\!Ax = Vb$) or as a right-preconditioner (solving $AVy = b$, $x = Vy$), potentially with different convergence characteristics, despite having identical spectral properties \cite{saad}. We tested both versions in this work and observed that the left- and right-preconditioned GMRES methods often resulted in nearly identical convergence rates; consequently, in our numerical experiments we focus just on the left-preconditioned case.\footnote{We prefer the left-preconditioned version as the corresponding residual used by GMRES is then a good proxy of the PDE (backward) error, since $V\!A$ usually has eigenvalues very close to 1.}

\subsection{Element-wise block smoothers}
\label{sec:smoothers}

We consider a variety of multigrid smoothers (also known as relaxation methods), specifically block Jacobi, block Gauss-Seidel, and block SAI smoothers. Here, each \textit{block} corresponds to the collective set of degrees of freedom on each mesh element. Accordingly, we partition $Ax = b$ into blocks $A_{ij}$ so that $\sum_j A_{ij} x_j = b_i$, where $(\cdot)_i$ denotes the set of coefficients (in our case, tensor-product Legendre polynomial coefficients) on element $i$ and $A_{ij}$ denotes the corresponding $(i,j)$-th block of $A$. To define each of the smoothers we consider just one application or iteration thereof, denoting the input approximation by $\xin \in V_h$ and the output (i.e., updated/smoothed) approximation by $\xout \in V_h$.

\subsubsection{Block Jacobi} One iteration of (damped) block Jacobi corresponds to the update
\begin{equation} \label{eq:Jsmooth} \xout_i := (1 - \omega)\, \xin_i + \omega\, A_{ii}^{-1} \bigl(b_i - \textstyle \sum_{j \neq i} A_{ij} \xin_j \bigr), \end{equation}
for every element $i$. For elliptic interface problems, each and every diagonal block $A_{ii}$ is a symmetric positive (semi)definite matrix; in our particular implementation, we precompute and cache a symmetric eigendecomposition of $A_{ii}$ and use this factorisation to perform the inversion in \cref{eq:Jsmooth}. Meanwhile, $\omega > 0$ is a damping parameter required for convergence---for DG methods especially, undamped Jacobi (corresponding to $\omega = 1$) usually results in a divergent method, unsuitable for multigrid. For the numerical study of this work, we have set $\omega = 0.8$, empirically determined to yield near-optimal convergence rates while also being somewhat conservative (in case, e.g., the Jacobi smoother is applied to highly unstructured meshes for which extra damping is sometimes needed).

\subsubsection{Block Gauss-Seidel} \label{sec:GS} Gauss-Seidel is algorithmically equivalent to Jacobi, except the updates are overwritten immediately in-place as the elements are visited in some particular order. Mathematically,
\[ \xout_i := (1 - \omega)\, \xin_i + \omega\, A_{ii}^{-1} \bigl(b_i - \textstyle \sum_{j \neq i} A_{ij} \mathring{x}_{j,i} \bigr) \quad \text{where} \quad \mathring{x}_{j,i} := \begin{cases} \xin_j & \text{if } i < j, \\ \xout_j & \text{if } j < i, \end{cases}  \]
where ``$i < j$'' means element $i$ is updated before element $j$. We consider two approaches:
\begin{enumerate}
    \item \textit{Multi-coloured Gauss-Seidel.} In a setup phase, a graph-colouring algorithm is applied to the connectivity graph associated with the block sparsity pattern of $A$ to define a colour for each element. (On Cartesian grids with one-sided numerical fluxes, this recovers the well-known red-black colouring associated with the 5-point (2D) or 7-point (3D) Laplacian stencil.) In the Gauss-Seidel sweep, every element of the same colour is processed (possibly in parallel), before moving onto the next colour. Multi-coloured Gauss-Seidel does not require any damping, so we use $\omega = 1$. 
    \item \textit{Processor-block Gauss-Seidel.} This method arises from a standard domain decomposition approach for distributed memory parallel computation (e.g., using MPI), see, e.g., \cite{AdamsBrezinaHuTuminaro}. Each subdomain is assigned to a processor, with a ghost layer of elemental values synchronised and then frozen at the beginning of each iteration. Within each subdomain, each processor applies standard Gauss-Seidel to update only its locally-owned elements, according to some predefined ordering. Processor-block Gauss-Seidel requires some amount of damping, which here we set to $\omega = 0.8$. (Indeed, in the limit of one element per processor, the method decays to block Jacobi.)
\end{enumerate}
We have included both Gauss-Seidel methods in our numerical study as they represent two of the most common implementations. Multi-coloured Gauss-Seidel yields the fastest convergence and (once coloured) allows for higher parallelism, e.g., in a multithreaded setting; on the other hand, the colouring process itself is in general a global operation whose implementation may be nontrivial. Processor-block Gauss-Seidel (possibly with a multi-coloured approach on each processor's subdomain) might be simpler to implement; on the other hand, it needs some damping to ensure reliable convergence, which may or may not offset its simplified inter-process communication patterns.

\subsubsection{Block SAI} \label{sec:SAI} Sparse approximate inverse methods solve a minimisation problem to find the best possible inverse with a given sparsity structure. One could also apply adaptive algorithms to iteratively refine the sparsity pattern, usually by some kind of controlled fill-in, so as to improve the approximation (e.g., to automatically handle local anisotropy) \cite{1997doi:10.1137/S1064827594276552,2000doi:10.1137/S106482759833913X,2000doi:10.1137/S0895479899339342,2001doi:10.1137/S1064827500380623,2002BROKER200261,2010doi.org/10.1155/2010/930218}. 
In this work, we focus on a straightforward method that uses a fixed sparsity structure matching the block-sparse pattern of the powers of $A$.\footnote{Although simple, this approach may nevertheless serve as a guide toward evaluating the potential benefits of a more complex/adaptive approach.} Specifically, let $\ell$ be a nonnegative integer and let ${\mathcal M}_{A}^\ell$ denote the set of matrices whose block-sparse pattern coincides with $A^\ell$. For example, ${\mathcal M}_{A}^0$ denotes the set of block-diagonal matrices, ${\mathcal M}_{A}^1$ denotes the set of matrices with the same block-sparse pattern as $A$, etc. This setup can also be viewed in terms of the elemental stencils associated with the discrete Laplacian operator: $\ell = 0$ corresponds to initialising with a trivial stencil (one block per element), $\ell = 1$ yields the same stencil as the discrete Laplacian, $\ell = 2$ arises from applying the Laplacian's stencil to itself thereby yielding a larger stencil, etc. Applied to $A$ unmodified, the prototype SAI approach corresponds to computing the SAI $B_\ell$ such that
\[ B_\ell := \argmin_{B \in {\mathcal M}_{A}^\ell} \| A B - \mathbb I \|_F, \]
where $\| \cdot \|_F$ denotes the Frobenius norm and $\mathbb I$ denotes the identity matrix. Note that every block column of $B_\ell$ can be computed independently of (and possibly concurrently with) every other block column. In particular, computation of the $i$-th block column of $B_\ell$ amounts to solving for its nonzero blocks through a relatively small, overdetermined, uniformly-weighted, least squares problem; in our implementation, we have used a simple dense least squares QR method from \textsc{lapack}. In fact, it is often the case many block columns of $B_\ell$ have the \textit{same} block-sparse solution, up to a one-to-one mapping of the blocks involved---each representative solution could therefore be computed once, cached, and reused across all appropriate columns of $B_\ell$, and even across the multigrid hierarchy itself, so as to speed up the cost of computing $B_\ell$ across all grid levels; more comments on this aspect are provided in the concluding remarks.

The above process defines an SAI directly from $A$, however it can sometimes be useful to precondition or ``balance'' the system. In particular, elliptic interface problems of the form \cref{eq:scalar1,eq:scalar2,eq:scalar3} often involve viscosity coefficients varying over multiple orders of magnitude. Accordingly, the various blocks $A_{ij}$ can have considerably different magnitudes. Since an ordinary SAI method solves a uniformly-weighted least squares problem, its solution may therefore be skewed leading to poor smoothing behaviour, essentially because nearby elements are smoothed at different rates. Fortunately, a simple method can be implemented that robustly handles large viscosity ratios and restores good smoothing behaviour. Our ``balancing'' approach performs a simple (temporary) diagonal pre- and post-scaling of $A$ so as to equalise the Frobenius norm of its diagonal blocks. Let $\alpha$ be a vector with the same dimension as $x$ and imagine it to have the same block structure; let $\alpha_{i,j}$ denote the $j$-th degree of freedom on element $i$. We define the values of $\alpha$ by
\[ \alpha_{i,j} := \| A_{ii} \|_F^{-1/2} \]
for every $i,j$, so that every diagonal block of $\tilde A := \diag(\alpha) A \diag(\alpha)$ has unit Frobenius norm. We then apply the standard SAI algorithm to $\tilde A$ and call\footnote{We sometimes drop the subscript $\ell$ whenever it is not crucial to specify.} the result $\tilde B$; the true SAI of $A$ is then defined by $B := \diag(\alpha) \tilde B \diag(\alpha)$. One can view this process as a simple diagonal preconditioning; alternatively, one can view it as modifying the standard SAI approach to use a weighted least squares procedure. Indeed, $\tilde A \tilde B \approx \mathbb I$, where the approximation denotes a minimisation in the unweighted Frobenius norm; on the other hand, $A B \approx \mathbb I$, but now the approximation involves a weighted Frobenius norm, controlled by $\alpha$, that effectively ``cancels'' the effect of large viscosity coefficients. Note that this balancing approach preserves the potential benefits of SAI; e.g., each block column of $B$ can still be (pre)computed independently of (and concurrently with) all others. We have used this balancing approach throughout the numerical results of this paper.

Once computed, an SAI $B$ can then be used as a multigrid smoother through one of two distinct forms: either through the update
\begin{equation} \label{eq:B1} \xout := \xin - B(A \xin - b), \end{equation}
or the update
\begin{equation} \label{eq:BT} \xout := \xin - B^\trans (A \xin - b). \end{equation}
Even when $A$ is symmetric, an SAI is generally not expected to be symmetric. Consequently, there can be a difference between using \cref{eq:B1} or \cref{eq:BT} as a multigrid smoother. In fact, our results show that some combinations lead to poor multigrid performance, while other combinations yield excellent results.

\subsection{Numerical experiments}
\label{sec:scalarexp}

In this section we examine multigrid performance for a variety of smoother and solver combinations, on a variety of elliptic interface problems, for a variety of polynomial degrees $\polydeg$, in two and three dimensions, across multiple grid sizes. Our primary metric used to assess multigrid efficiency concerns the convergence rate of multigrid-preconditioned CG or GMRES. Specifically, we randomly generate a right hand side vector $b$ and solve $Ax = b$ using a zero initial guess $x_0 \equiv 0$. Each iteration $k$ of multigrid-preconditioned CG or GMRES yields an improved solution $x_k$ with corresponding residual norm $r_k$; typically, $r_k$ is a nearly-constant fraction of $r_{k-1}$ so that $r_k \approx \rho\, r_{k - 1}$, where $\rho$ is the convergence rate. Using this empirically determined convergence rate,\footnote{The coefficients of $b$ are randomly and independently drawn from the uniform distribution on $[-1,1]$. With high probability the corresponding problem $Ax = b$ contains modes for which convergence is the slowest, thereby estimating worst-case convergence rates. The residual norm is the same as that used by the particular Krylov method: in the case of left-preconditioned GMRES, the residual is equal to $r_k := \|VAx_k - Vb\|_2$, while for preconditioned CG it is computed via the $V$ inner product, $r_k^2 := (Ax_k - b)^\trans V (Ax_k - b)$. On a log-linear graph, $r_k$ is very well approximated by a line of negative slope; we calculate the value of $\rho$ via simple linear regression on this graph. Our numerical tests typically have 5 to 10 data points, except when convergence is so fast that just a few points are collected before reaching the limits of double-precision arithmetic.} our primary metric used to assess multigrid speed is defined by
\begin{equation} \label{eq:eta} \eta := \frac{\log 0.1}{\log \rho}. \end{equation}
Here, $\eta$ represents the expected number of iterations needed to reduce the residual by a factor of 10; e.g., $k \eta$ represents the number of iterations to reduce it by $10^k$. We prefer to use $\eta$ instead of $\rho$, as the former is intuitively more useful when comparing different solvers; e.g., a two-fold reduction in $\eta$ coincides with a two-fold reduction in iteration count and thus, all else being equal, two times faster.

We test the performance of a variety of multigrid solver combinations. In each case, the V-cycle of \cref{algo:vcycle} is used as the preconditioner for CG or as a left-preconditioner for GMRES, using $\nu_1 = \nu_2 = 3$ pre- and post-smoothing steps. We consider the following smoother options:
\begin{itemize}
    \item \textit{Block Jacobi.} Applying (damped) block Jacobi as the pre- and post-smoother automatically results in a symmetric V-cycle. Combined with CG and GMRES, we denote the resulting solver by \xJCG{} and \xJGL, respectively.

    \item \textit{Block Gauss-Seidel.} By itself, Gauss-Seidel is not symmetric, owing to its dependence on the ordering of elemental updates. As such, there are two basic approaches for using Gauss-Seidel in a V-cycle: (i) use the \textit{same} predefined, fixed ordering for \textit{both} the pre-smoothing step and the post-smoothing step; or (ii) use a predefined, fixed ordering in the pre-smoothing step, and then the \textit{reverse} ordering in the post-smoothing step. Option (ii), denoted by \xOI, leads to a symmetric V-cycle and can be used with both CG and GMRES; in contrast, option (i), denoted by \xOO, results in an asymmetric V-cycle and can only be used with GMRES. Besides these combinations, we also have the choice of multi-coloured Gauss-Seidel, denoted by \xGS, and processor-block Gauss-Seidel, denoted by \xLGS, as described in \cref{sec:GS}. In total this leads to six combinations, denoted by \xGSOICG{} and \xLGSOICG{} when coupled to CG, and \xGSOIGL{}, \xGSOOGL{}, \xLGSOIGL{}, and \xLGSOOGL{} when coupled to GMRES.

    \item \textit{Block SAI.} We consider three SAI solvers, \SAI0, \SAI1, and \SAI2, where the integers correspond to the value $\ell$ defining the block-sparse pattern of the approximate inverse, see \cref{sec:SAI}. Analogous to the Gauss-Seidel case, there are two basic approaches for using SAI in a V-cycle: (i) use either \cref{eq:B1} or \cref{eq:BT} for the pre-smoothing step and the \textit{same} choice for the post-smoothing step; or (ii) use one option for the pre-smoothing step and the \textit{opposite} option for the post-smoothing step. Option (ii), denoted by \xOI{} or \xIO{}, results in a symmetric V-cycle and can be used with both CG and GMRES; in contrast, option (i), denoted by \xOO{} or \xII, results in an asymmetric V-cycle and can only be used with GMRES. The arrow symbols used here are in analogy with Gauss-Seidel sweeps: applying \cref{eq:B1} is in some sense a ``forward'' sweep, whereas applying \cref{eq:BT} ``reverses'' the direction; e.g., \xOI{} means to apply \cref{eq:B1} as pre-smoother and \cref{eq:BT} as post-smoother. In total, these various solver combinations lead to six possibilities for each $\ell$, denoted by \xSaiOICG{} and \xSaiIOCG{} when coupled to CG, and \xSaiOIGL{}, \xSaiIOGL{}, \xSaiOOGL{}, and \xSaiIIGL{}, when coupled to GMRES.
\end{itemize}
In the presented results, we will selectively prune these combinations as we progress from simpler to more complex cases.

\begin{table}%
\centering%
\sffamily\small%
\newcommand{\incsquare}[1]{\raisebox{-0.1\height}{\setlength{\fboxsep}{0pt}\setlength{\fboxrule}{0.1pt}\fbox{\includegraphics[width=1em]{#1}}}}
\newcommand{\bul}[1]{\raisebox{#1 ex}{\scalebox{0.7}{$\bullet$}}}
\begin{tabular}{c@{\hspace{0.5mm}}c@{\hspace{1mm}}l@{\hspace{10mm}}c@{\hspace{1mm}}l@{\hspace{10mm}}c@{\hspace{1mm}}l}
\xCG & \xGMRES & {\footnotesize Solver/smoother} \\
\incsquare{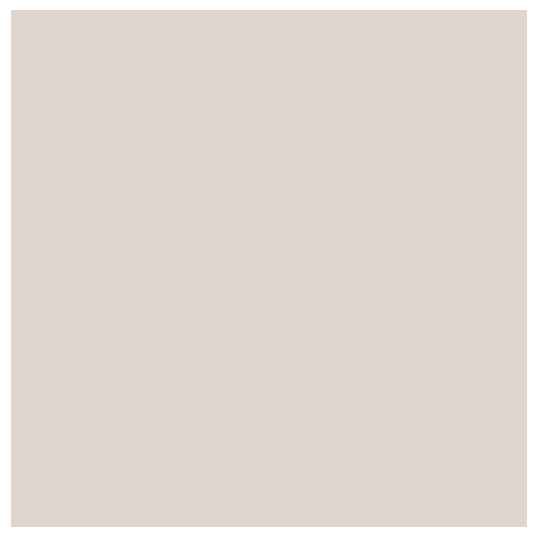} & \incsquare{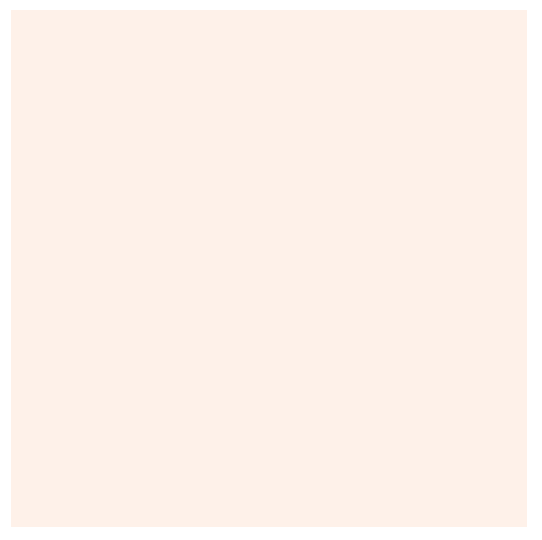} & Jacobi &
\raisebox{-0.2ex}{$\xOO \xOI$} & \multirow{2}{*}{Pre/post ordering} &
$\bullet$ & $\eta \in (0,2)$  \\
\incsquare{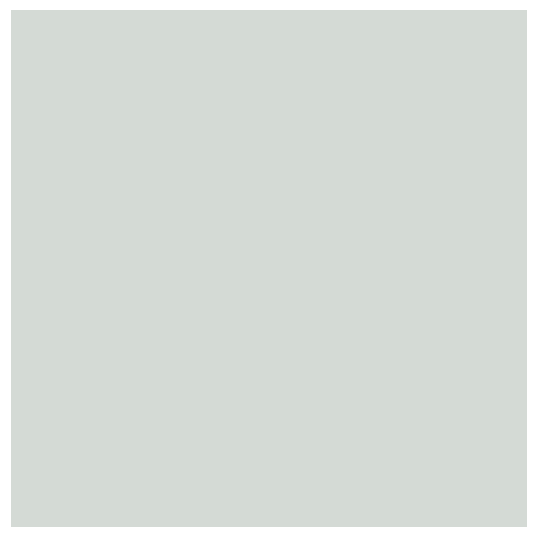} & \incsquare{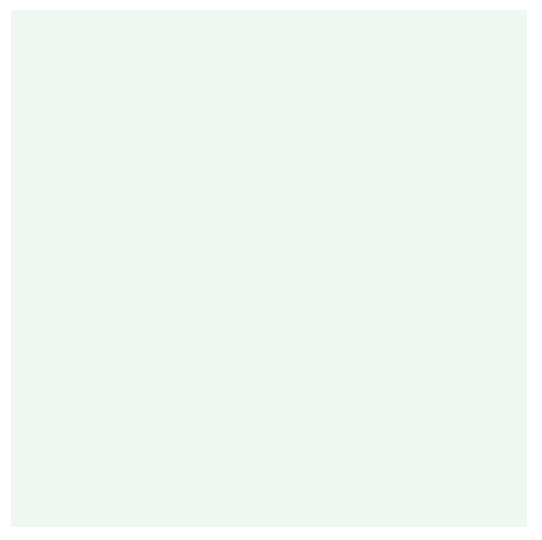} & \xLGSsymbol~Gauss-Seidel &
\raisebox{0.2ex}{$\xIO \xII$} & & 
$\scalebox{1}[0.85]{$\blacktriangle$}$ & $\eta \in [2,10)$ \\
\incsquare{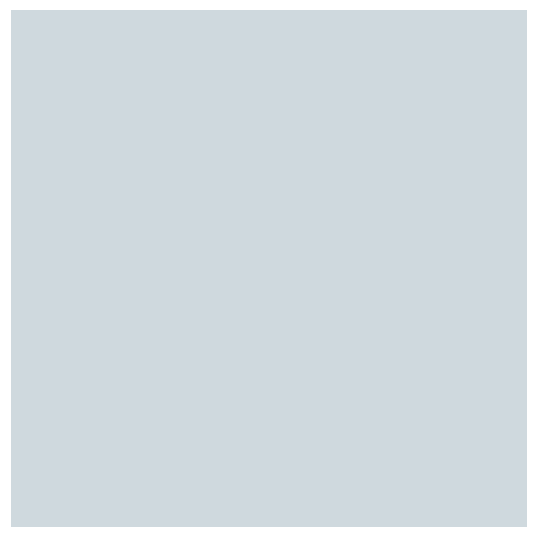} & \incsquare{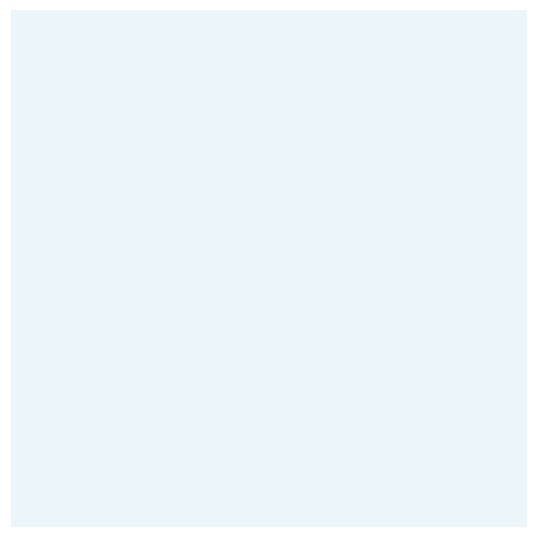} & \xGSsymbol~Gauss-Seidel &
\bul{0.2}\bul{0.4}\bul{0.6}\bul{0.6}\bul{0.4}\bul{0.5} & $n = 8,16,32,64,\ldots$ &
$\times$ & $\eta \in [10,\infty)$ \\
\incsquare{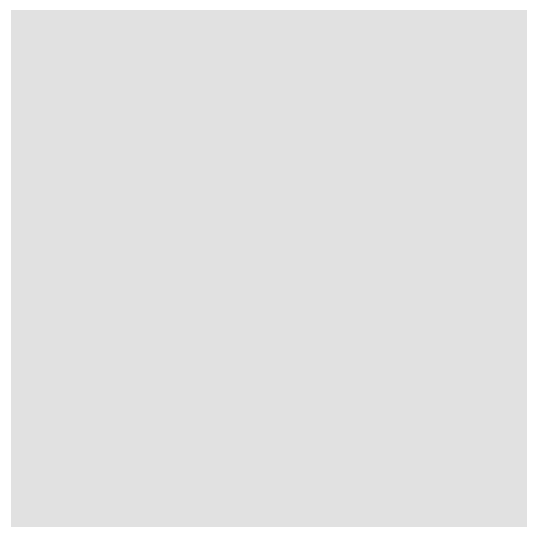} & \incsquare{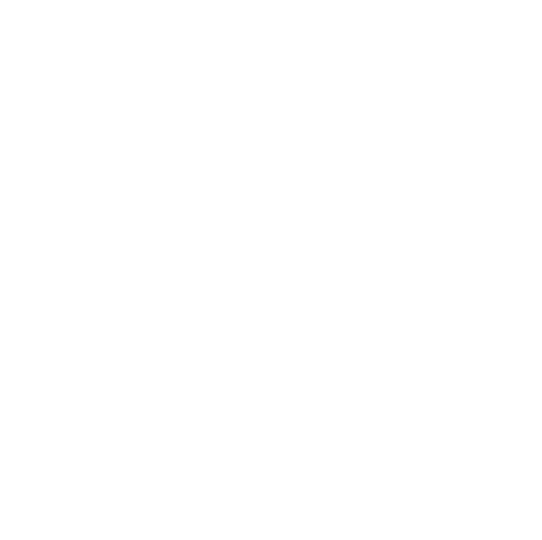} & $\SAI{\ell}$
\end{tabular}
\captionsetup{name=Legend}
\caption{Legend for the numerical results of \cref{fig:01} through \cref{fig:11}. Darker shaded (resp., lighter or unshaded) columns correspond to multigrid-preconditioned CG (resp., GMRES). Arrows indicate the choice of ordering in pre-smoothing and post-smoothing steps, as defined in the in the paragraph following \Cref{eq:eta}. Each group of symbols plot the numerically determined speed $\eta$ on grid sizes $n \times n\,(\times\, n)$ where $n = 8, 16, 32, \ldots$ from left-to-right.} 
\label{tab:legend}
\end{table}

\newcommand{\seelegend}{A visualisation guide is provided in \hyperref[tab:legend]{Legend 1} on page \pageref{tab:legend}.}

For each specific PDE problem, across the various solver combinations, we also consider a variety of polynomial degrees $\polydeg$, across multiple grid sizes. Ideally, the number of CG/GMRES iterations required to reach some prescribed residual reduction threshold should stabilise and remain roughly constant as the mesh size is progressively refined. In 2D, grid sizes range from $8 \times 8$ to $1024 \times 1024$, with polynomial degrees $\polydeg \in \{1, 2, \ldots, 5\}$; in 3D, we consider $\polydeg \in \{1, 2, 3\}$ and grid sizes up to $256 \times 256 \times 256$. To convey this abundance of information, our numerical results are presented through graphs of the form shown in \cref{fig:01}. Each row corresponds to a specific spatial dimension and polynomial degree, while each vertical column corresponds to a specific solver/smoother combination, as indicated at the bottom of the figure; in the set of columns associated with the same smoother type, those with a darker shading correspond to a CG solver, whereas lighter or unshaded columns correspond to GMRES. Each individual group of dots plots the numerically determined multigrid speed $\eta$ across multiple grid sizes: from left-to-right, each dot plots $\eta$ for grids of size $n \times n\,(\times\, n)$ with $n = 8, 16, 32, \ldots$ in that order. Sometimes a particular solver performs poorly and results in slow convergence: we consider multigrid speeds with $2 \leq \eta < 10$ to be ``slow'' and when this occurs, the corresponding dot is replaced by a triangle $\scalebox{1}[0.85]{$\blacktriangle$}$ at the top of the graph; if $\eta \geq 10$ (corresponding to $\rho \gtrapprox 0.8$) we declare the method to be nonconvergent (at least, from a practical point of view) and indicate this with a cross $\times$. For reference, \hyperref[tab:legend]{Legend 1} provides a summary of this visualisation scheme.

\subsubsection{Classical Poisson problem}
\label{sec:poisson}

\newcommand\orbscale{0.806}
\begin{figure}%
\centering%
\includegraphics[scale=\orbscale]{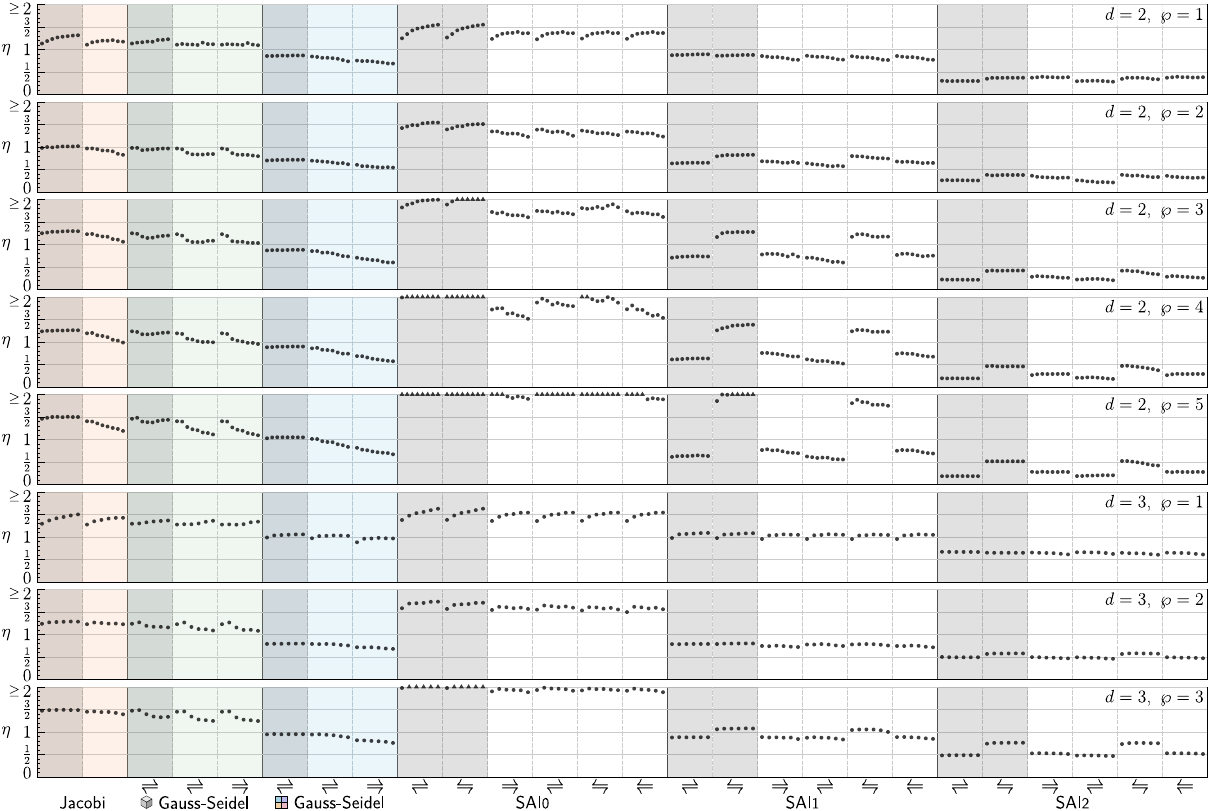}%
\caption{Multigrid solver performance for the classical single-phase Poisson problem considered in \cref{sec:poisson}. \seelegend{}}%
\label{fig:01}%
\end{figure}

We begin with perhaps the simplest case, that of a (single-phase) Poisson problem on the unit rectangular domain $\Omega = (0,1)^d$ with $\mu \equiv 1$, applying periodic boundary conditions and uniform Cartesian grids. \cref{fig:01} contains the corresponding numerical results for all of the tested solver combinations. Several patterns emerge, roughly representative of all the results presented in this work: (i) unsurprisingly, when using \SAI{\ell} as smoother, multigrid convergence rates get progressively faster as $\ell$ increases, owing to the approximate inverse becoming more and more accurate; (ii) among all smoothers, \SAI0 is usually the slowest and worsens as $\polydeg$ increases; (iii) for a fixed $\polydeg$, out of the two block Gauss-Seidel approaches, the multi-coloured (undamped) method is fastest, as expected, while the speed of processor-block Gauss-Seidel sits between that of $\xGS$ and (damped) block Jacobi, depending on the parallel runtime;\footnote{The numerical results in this work have been collected using a hybrid MPI/multithreading framework. 2D problems are sometimes executed without MPI, sometimes with MPI on a few or perhaps hundreds of ranks; 3D problems sometimes run on a few hundred ranks, other times on tens of thousands. Since processor-block Gauss-Seidel smoothing behaviour depends on the subdomain layout, this means there can be some variance in its measured convergence rates $\eta$.} (iv) the best \SAI1 method has speed comparable to the best Gauss-Seidel method; (v) across all solver combinations, \SAI2 generally exhibits the fastest speeds $\eta$; and (vi) besides the standout cases where solver speed is noticeably slow, every other solver combination demonstrates classic ``textbook'' multigrid behaviour whereby the CG/GMRES iteration count remains roughly constant as the grid is refined. Beyond these broad observations there are some more specific ones:
\begin{itemize}
    \item Among the \xGS{} solvers and the corresponding pre-/post-smoother combinations, \xGSOOGL{} is usually the fastest. This is an example where, although a good symmetric V-cycle preconditioned CG solver is available to us (with comparatively smaller memory footprint), it can be noticeably faster (e.g., two thirds as many iterations in some cases) to use an asymmetric V-cycle preconditioned GMRES solver. Intuitively, we interpret this speedup as follows. Consider \xGS\xOI{} in the case of a simple red-black colouring: the pre-smoother visits the elements in the order \textsf{RBRBRB} followed by the post-smoother in the order \textsf{BRBRBR}; in between is a coarse grid correction, but if that correction is negligible (resp., zero), then there are two black sweeps in a row with the second sweep having negligible effect (resp., no effect whatsoever). Since \xGS\xOO{}, which keeps the alternating \textsf{RB} pattern from pre- to post-smoothing, is generally observed to be faster than \xGS\xOI{}, we view the latter as ``wasting'' some effort in smoothing the same set of elements immediately before and after the coarse grid correction.
    \item Among the \SAI1 and \SAI2 solvers, \xSaiIOCG{} and \xSaiIOGL{} clearly standout as the worst performers. In other words, it is suboptimal to use \cref{eq:BT} as a pre-smoother and \cref{eq:B1} as a post-smoother.
\end{itemize}
In addition to periodic boundary conditions, we also tested Dirichlet and Neumann boundary conditions. The results are quite similar to those indicated in \cref{fig:01} and thus not included here. Essentially the only difference is that, in the case of Neumann boundary conditions, the above noted bad performance of \xSaiIOCG{} and \xSaiIOGL{} is amplified, particularly so for $\ell = 1$.

\subsubsection{Adaptive mesh refinement} 
\label{sec:amr}

\begin{figure}%
\centering%
\raisebox{-0.5\height}{\includegraphics[height=1.5in]{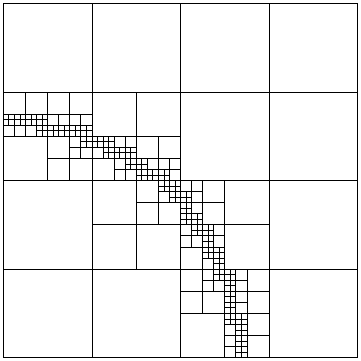}} \qquad \raisebox{-0.5\height}{\includegraphics[height=1.75in]{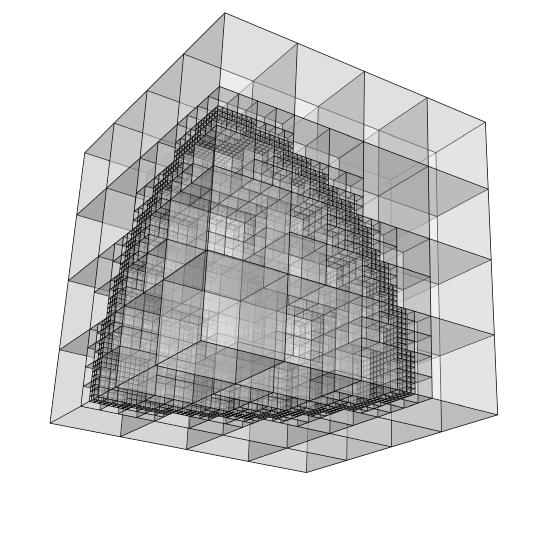}}%
\caption{Examples of the nonuniform quadtree (left) and octree (right) used in the adaptive mesh refinement problem in \cref{sec:amr}. Each example corresponds to a primary grid size of $n = 4$, refined up to four levels around a circular arc or spherical shell of the same radius.}%
\label{fig:amr}%
\end{figure}

\begin{figure}%
\centering%
\includegraphics[scale=\orbscale]{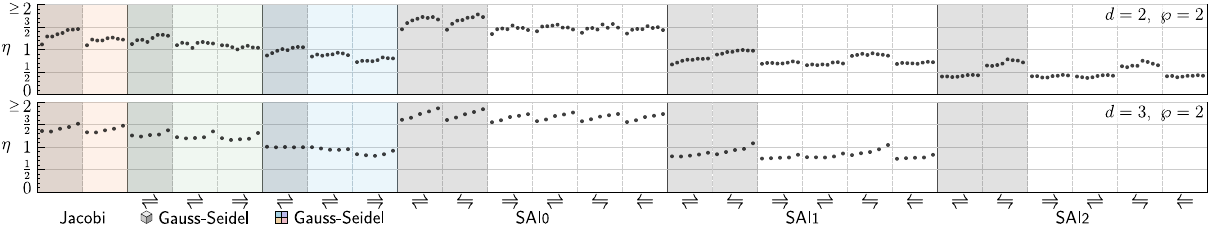}%
\caption{Multigrid solver performance for the quadtree/octree nonuniform mesh problem considered in \cref{sec:amr}. \seelegend{}}%
\label{fig:04}%
\end{figure}

Our second test problem is similar to the first, except we now apply it to a nonuniform adaptively refined mesh that has up to four levels of nested refinement around a circular arc (in 2D) or spherical shell of the same radius (in 3D), as illustrated in \cref{fig:amr}. On this mesh we consider a (single-phase) Poisson problem with $\mu \equiv 1$ and Neumann boundary conditions. Owing to the somewhat rapid refinement, the block sparsity pattern of $A^2$ is substantially more dense compared to the case of a uniform grid; this makes the $\ell = 2$ SAI appreciably more expensive to build (mainly due to the cubic scaling in the number of unknowns when solving the associated least squares problems) --- on a high resolution 3D mesh this unfortunately renders the \SAI2 solver too computationally costly, so we exclude \SAI2 in the 3D case. \cref{fig:04} compiles the results, measuring multigrid solver speed across a variety of solver choices, focusing on the case of $\polydeg = 2$ for simplicity. Overall, compared to the case of uniform grids, we observe similar results; in particular, the best \SAI1 solver has about 10--20\% fewer iterations than the best Gauss-Seidel solver.

\subsubsection{Multiphase elliptic interface problems}
\label{sec:multiphase}

\begin{figure}%
\centering%
\includegraphics[width=\textwidth]{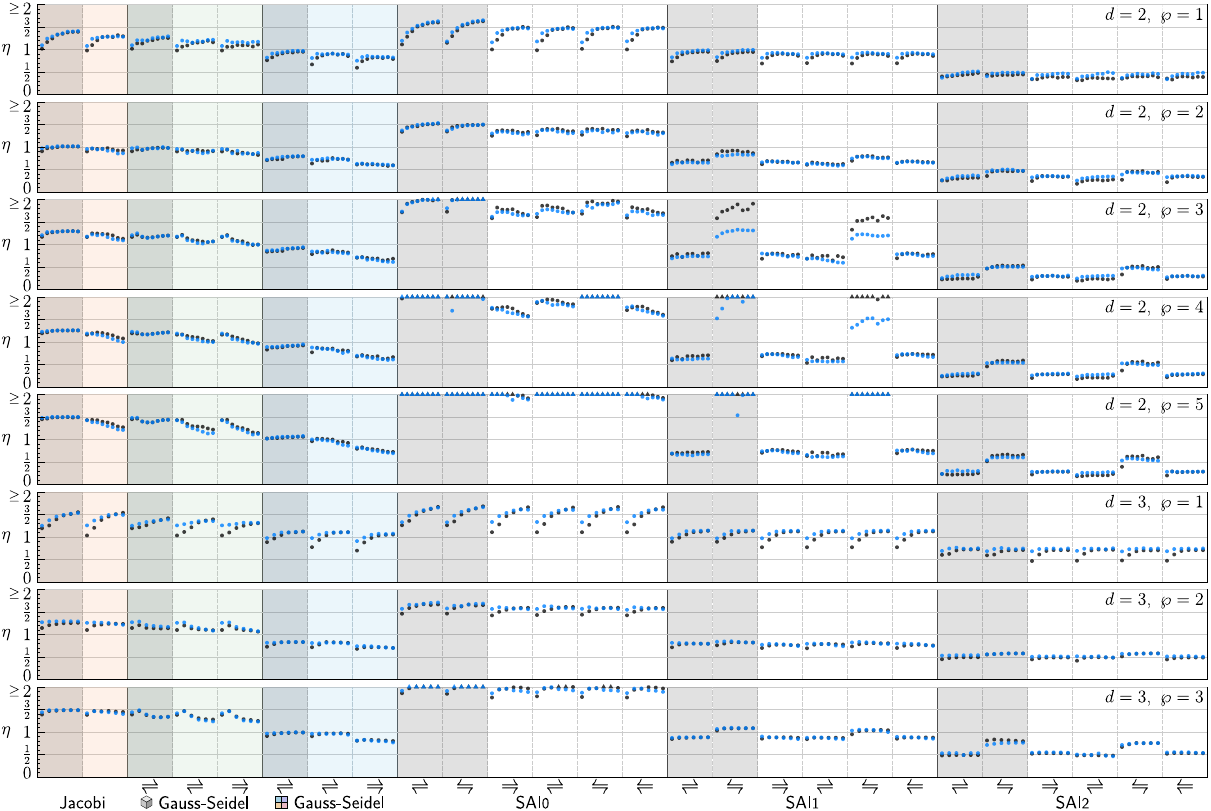}%
\caption{Multigrid solver performance for the multiphase elliptic interface problem considered in \cref{sec:multiphase}. \seelegend{} Black (resp., blue) symbols correspond to the viscosity ratio $\mu_0/\mu_1 = 10^{-4}$ (resp., $10^4$). Other viscosity ratios, including $10^{\pm2}$, $10^{\pm6}$ and $10^{\pm8}$, yield nearly identical results.}%
\label{fig:05}%
\end{figure}

We now consider a potentially more challenging case of a two-phase elliptic interface problems wherein the ellipticity coefficient $\mu$ exhibits a jump several orders in magnitude across an embedded interface. Let $\Omega = (0,1)^d$ be divided into an interior rectangular phase $\Omega_0 = (\tfrac14, \tfrac34)^d$ and an exterior phase $\Omega_1 = \Omega \setminus \overline{\Omega_0}$, each with viscosity $\mu_i$. Two different viscosity ratios are evaluated: $\mu_0/\mu_1 = 10^4$ and $\mu_0/\mu_1 = 10^{-4}$. We consider uniform Cartesian grids and periodic boundary conditions on $\partial \Omega$ (similar results are obtained for Dirichlet or Neumann boundary conditions). \cref{fig:05} compiles the multigrid speed results for the two viscosity ratios, for all of the same parameters and combinations as was considered in our first test case. As before, with the exception of many \SAI0 cases and the outliers \xSaiIOCG{} and \xSaiIOGL{}, we observe asymptotically stable multigrid solver speeds $\eta$. This favourable outcome is attributed in large part to the use of viscosity-upwinded numerical fluxes \cite{fluxx}. In particular, the cited work explored multigrid performance mainly in the context of Gauss-Seidel methods; the results of \cref{fig:05} show that viscosity-upwinded fluxes also work effectively in combination with damped block Jacobi as well as \SAI1 and \SAI2 smoothers. Finally, we note that other viscosity ratios were also tested, including $10^{\pm2}$, $10^{\pm6}$ and $10^{\pm8}$; across all of these the results are essentially identical to those shown in \cref{fig:05}.

\section{Multiphase Stokes Problems}

Thus far, we have explored efficient multigrid solvers for LDG discretisations of scalar elliptic interface problems and showed that various choices yield the ``gold standard'' performance of classical geometric multigrid methods. In this section, we consider the arguably more challenging case of multiphase Stokes problems and demonstrate that a similar outcome can be achieved.

Our prototype (steady-state) Stokes problem consists of solving for a velocity field $\vu : \Omega \to \R^d$ and a pressure field $p: \Omega \to \R$ such that 
\begin{equation} \label{eq:govern1} \left. \begin{aligned} -\nabla \cdot \bigl(\mu_i (\nabla \vu + \gamma \nabla \vu^\trans) \bigr) + \nabla p &= \mathbf f \\ -\nabla \cdot \mathbf \vu &= f \end{aligned} \right\} \text{ in } \Omega_i, \end{equation}
subject to the interfacial jump conditions
\begin{equation} \label{eq:govern2} \left. \begin{aligned} \jump{\vu} &= {\mathbf g}_{ij} \\ \jump{\mu (\nabla \vu + \gamma \nabla \vu^\trans) \cdot \vn - p \vn} &= {\mathbf h}_{ij} \end{aligned} \right\} \text{ on } \Gamma_{ij}, \end{equation}
and boundary conditions
\begin{equation} \label{eq:govern3} \begin{aligned} \vu &= {\mathbf g}_{\partial} && \text{on } \Gamma_D, \\ \mu (\nabla \vu + \gamma \nabla \vu^\trans) \cdot \vn - p \vn &= {\mathbf h}_{\partial} && \text{on } \Gamma_N, \end{aligned} \end{equation}
where $\Omega$ is a domain in $\mathbb R^d$ divided into one or more subdomains/phases $\Omega_i$, $\Gamma_{ij} := \partial \Omega_i \cap \partial \Omega_j$ is the interface between phase $i$ and $j$, and $\Gamma_D$ and $\Gamma_N$ denote the parts of $\partial \Omega$ upon which velocity Dirichlet or stress boundary conditions are imposed, respectively. Here, $\gamma \in \{0,1\}$ is a parameter defining the form of the Stokes equations: they are said to be in \textit{standard form} when $\gamma = 0$ and \textit{stress form} when $\gamma = 1$. Finally, $\mu_i$ is a phase-dependent viscosity coefficient, while $\vf$, $f$, $\vg$, and $\vh$ provide the data for the multiphase Stokes problem and are given functions defined on $\Omega$, its boundary, and internal interfaces.

\subsection{LDG schemes}

As before, we consider LDG discretisations on Cartesian grids, quadtrees, or octrees, along with tensor-product polynomial spaces. Solving for a discrete solution $\vu_h \in V_h^d$ and $p_h \in V_h$, the high-order LDG method for \cref{eq:govern1,eq:govern2,eq:govern3} is built according to a five-step process: (i) using viscosity-upwinded one-sided numerical fluxes, define a discrete stress tensor $\vtau_h \in V_h^{d \times d}$ equal to the discretisation of $\nabla \vu_h + \gamma \nabla \vu_h^\trans$, taking into account Dirichlet source data $\vg$; (ii) define $\vsigma_h \in V_h^{d \times d}$ as the viscous stress $\mu \vtau_h - p_h \mathbb I$ via an $L^2$ projection of $\mu \vtau_h$ onto $V_h^{d \times d}$; (iii) using the adjoint of the discrete gradient operator from the first step, compute the negative discrete divergence of $\vsigma_h$, taking into account Neumann-like data $\vh$, and set the result equal to the $L^2$ projection of $\vf$ onto $V_h^d$; (iv) compute the discrete divergence of $\vu_h$, taking into account Dirichlet data $\vg \cdot \vn$, and set the result equal to the $L^2$ projection of $f$ onto $V_h$; last, (v) amend the overall system by appropriately including penalty stabilisation. Our LDG framework is based upon the scheme initially developed by Cockburn \textit{et al} \cite{CockburnKanschatSchotzauSchwab2002}, in part extending it to the unsteady case as well as to multiphase problems. Full details of this construction are given in \cite{flame} as well as convergence studies demonstrating high-order accuracy (optimal order $\polydeg + 1$ in velocity and at least order $\polydeg$ in pressure, in the maximum norm). To summarise, the LDG scheme results in a symmetric linear system of the following form (here illustrated for the 3D case, with that of the 2D case being an obvious reduction),
\begin{equation} \label{eq:components}
\left( \begin{array}{@{}c@{\hspace{2mm}}c@{\hspace{2mm}}c|c@{}} \mathscr{A}_{11} & \mathscr{A}_{12} & \mathscr{A}_{13} & -G_1^\trans M \\ \mathscr{A}_{21} & \mathscr{A}_{22} & \mathscr{A}_{23} & -G_2^\trans M \\ \mathscr{A}_{31} & \mathscr{A}_{32} & \mathscr{A}_{33} & -G_3^\trans M \\ \hline -M G_1 & -M G_2 & -M G_3 & -E_p \end{array} \right)\!\!\begin{pmatrix} u_h \\ v_h \\ w_h \\ p_h \end{pmatrix} = \begin{pmatrix} b_u \\ b_v \\ b_w \\ b_p \end{pmatrix}\!, \end{equation}
where $u_h, v_h, w_h \in V_h$ denote the components of $\vu_h \in V_h^d$ while $(b_u,b_v,b_w,b_p)$ collects the entire influence of the source data $\mathbf f$, $f$, $\mathbf g$, and $\mathbf h$ onto the right hand side. Here, $M$ is a block-diagonal mass matrix, $G : V_h \to V_h^d$ is a discrete gradient operator, and $\mathscr{A} : V_h^d \to V_h^d$ implements the viscous part of the Stokes momentum equations, whose $(i,j)$-th block is given by 
\[ \mathscr{A}_{ij} = \delta_{ij} \bigl( \textstyle{\sum_{k=1}^d} G_k^\trans M_\mu G_k \bigr) + \gamma\, G_j^\trans M_\mu G_i + \delta_{ij} E, \]
where $M_\mu$ is a $\mu$-weighted block diagonal mass matrix and $\delta_{ij}$ the Kronecker delta. Here, $E$ is the same penalty operator encountered earlier for scalar elliptic problems; in contrast, $E_p$ is a penalty operator solely for the pressure variable and can play a crucial role in achieving fast multigrid performance, as discussed in detail in \cite{flame}. \Cref{eq:components} gives the system in component form; we will also refer to it through an abridged notation per the \{velocity, pressure\}-blocking of \cref{eq:components}, i.e.,
\begin{equation} \label{eq:blockform} \begin{pmatrix} \mathscr{A} & \mathscr{G} \\ -\mathscr{D} & -E_p \end{pmatrix} \begin{pmatrix} \vu_h \\ p_h \end{pmatrix} = \begin{pmatrix} {\mathbf b}_{\vu} \\ b_p \end{pmatrix}\!, \end{equation}
where $\mathscr G : V_h \to V_h^d$ is an effective discrete gradient operator and $\mathscr D : V_h^d \to V_h$ an effective discrete divergence operator. In terms of multigrid design, we shall also refer to it through the condensed form
\[ A_h x_h = b_h, \]
where $x_h$ collects $\vu_h$ and $p_h$ into one set of unknowns and $A_h$ is the entire symmetric saddle point operator. (The subscript $h$ may be dropped if context permits.)

\subsection{Multigrid schemes}

As for scalar elliptic interface problems, we consider multigrid methods using block Jacobi, block Gauss-Seidel, and block SAI smoothers. In the case of Stokes problems, we define each \textit{block} to correspond to the collective set of degrees of freedom on each mesh element, i.e., velocity and pressure combined. Accordingly, we partition $A x = b$ so that $\smash{\sum_j A_{ij} x_j = b_i}$, where $(\cdot)_i$ denotes the set of velocity and pressure values on element $i$ and $A_{ij}$ denotes the corresponding $(i,j)$-th block of $A$. With this setup, the block Jacobi and block Gauss-Seidel smoothers function equivalently to those detailed in \cref{sec:smoothers}. For example, block Gauss-Seidel sweeps over the elements, in some particular order, replacing $\smash{x_i \leftarrow A_{ii}^{-1} \bigl( b_i - \sum_{j \neq i} A_{ij} x_j \bigr)}$. Here, $A_{ii}$ is the $i$-th diagonal block of $A$ and takes on the form of a miniature Stokes operator; its inversion is a small saddle point problem whose wellposedness is generally guaranteed.\footnote{%
In general, a suitably-stabilised LDG Stokes method satisfies the inf-sup conditions \cite{CockburnKanschatSchotzauSchwab2002}.} In our particular implementation, we precompute and cache a symmetric eigendecomposition of $A_{ii}$ for every element $i$ and use this factorisation in the block Jacobi and block Gauss-Seidel sweeps.

\subsection{Balancing}
\label{sec:stokesbalance}

Recall that balancing is important for multiphase problems involving viscosities of different magnitudes. In fact, for Stokes problems specifically, balancing can be crucial even when the entire collection of elements have the same viscosity. To see why, note that the various operators within \cref{eq:blockform} scale with $h$ and $\mu$ in different ways: the viscous operator scales as $\mathscr{A} \sim \mu h^{d-2}$, pressure gradient operator scales as $\mathscr{G} \sim h^{d-1}$, velocity divergence operator scales as $\mathscr{D} \sim h^{d-1}$, and pressure penalty operator scales as $E_p \sim h^d/\mu$. Left untreated, sufficiently large $\mu$ would bias the least squares SAI solution to focus on the viscous operator block and essentially ignore the other blocks; conversely, sufficiently small $\mu$ would cause the SAI solution to focus too much on pressure stabilisation and ignore the other blocks. Note also the various blocks have different scalings in $h$, which means the length scales of the domain, or even the particular level in the multigrid hierarchy, also play a role and could adversely affect the SAI solution. Fortunately, all of these aspects can be resolved via a simple balancing process, similar to that used in \cref{sec:SAI}. In fact, one can tune this balancing to the benefit of overall multigrid speed effectively so that the SAI method smooths velocity and pressure at approximately the same speed.

In a form analogous to \cref{eq:blockform}, the $i$-th diagonal block of $A$ is equal to
\begin{equation} \label{eq:diagonalblock} \begin{pmatrix} \mathscr{A}_{ii} & \mathscr{G}_{ii} \\ -\mathscr{D}_{ii} & - E_{p,ii} \end{pmatrix}\!. \end{equation}
Our balancing approach performs a simple diagonal pre- and post-scaling of $A$ in order to control the norms of the four (newly-scaled) blocks in \cref{eq:diagonalblock}. Intuitively, we expect optimal results when these blocks have similar norms, so that the uniformly weighted least squares methods underlying SAI yield a smoother that equally balances the velocity and pressure degrees of freedom. Clearly, many choices exist for which norm to use and how precisely to perform this scaling. %
Satisfyingly, a relatively simple method was found that yielded optimal convergence rates in both 2D and 3D, on both standard- and stress-form Stokes problems, and on multiple polynomial degrees $\polydeg$. That method is to scale the system such that the top-left block in \cref{eq:diagonalblock} has unit Frobenius norm and such that the two anti-diagonal blocks have Frobenius norm about $0.1$. The value of $0.1$ is near optimal in most cases, i.e., roughly independent of $d$, $\gamma$, and $\polydeg$, but it can also be finely-tuned according to these parameters.

To precisely define this scaling, some additional notation is useful. Let $\alpha$ be the vector associated with pre- and post-diagonal scaling of $A$. We imagine it to have the same elemental block structure so that $\alpha_i$ denotes the collective set of degrees of freedom on element $i$. Define additional subscripting as follows: let $\alpha_{i,\vu,j,k}$ denote the component of $\alpha$ associated with the $k$-th degree of freedom of the $j$-th component of the velocity on the $i$-th element; similarly, let $\alpha_{i,p,k}$ denote the component of $\alpha$ associated with the $k$-th degree of freedom for the pressure on the $i$-th element.\footnote{Beyond the grouping into element-wise blocks, none of the algorithms described here depend on how these degrees of freedom are actually ordered.} We define the values of $\alpha$ as follows:
\begin{align*} \alpha_{i,\vu,j,k} &:= \|\mathscr{A}_{ii}\|_F^{-1/2},\\
\alpha_{i,p,j} &:= \zeta \|\mathscr{A}_{ii}\|_F^{1/2} \| \mathscr{D}_{ii}\|_F^{-1}, \end{align*}
for all $i, j, k$. Here, $\zeta \approx 0.1$ is a user-defined tunable value. It is straightforward to show that $\diag(\alpha) A \diag(\alpha)$ has the sought-after scaling: specifically, the top-left block of the (newly-scaled) system in \cref{eq:diagonalblock} has Frobenius norm one, while the two anti-diagonal blocks each have Frobenius norm $\zeta$. (This process also implicitly rescales the remaining bottom-right block associated with pressure stabilisation, though the outcome depends on the value of the corresponding penalty parameter; for the implementation choices in this work, the bottom-right block has a rescaled Frobenius norm typically between $0.2 \zeta$ and $0.7 \zeta$.)

We use $\alpha$ to build a balanced SAI smoother in the same way as discussed in \cref{sec:SAI}. That is, we build a proxy SAI, $\tilde{B}$, such that $\bigl[ \diag(\alpha) A \diag(\alpha) \bigr] \tilde{B} \approx \mathbb I$ (per the uniformly-weighted least squares procedure), and then define the true SAI of $A$ via $B = \diag(\alpha) \tilde{B} \diag(\alpha)$, so that $A B \approx \mathbb I$. Similar to the case of scalar elliptic interface problems, one can view this procedure as building an SAI using a weighted least squares procedure. The implementation only requires a small change so as to replace the balancing of each elemental block with the balancing of the $2 \times 2$ sub-blocks. In particular, the potential benefits of SAI are preserved: e.g., each block column of $B$ can be (pre)computed independently of (and concurrently with) all others, including also the possibility of caching these solutions for use across multiple columns and across the multigrid hierarchy.

\begin{table}%
\centering%
\sffamily\footnotesize%
\begin{tabular}{lcccccc}
& & \multicolumn{5}{c}{Polynomial degree $\polydeg$} \\
Stokes form & $d$ & 1 & 2 & 3 & 4 & 5 \\
\midrule
\multirow{2}{*}{Standard ($\gamma = 0$)} & 2D & 0.0382 & 0.131 & 0.0990 & 0.153 & 0.142 \\ 
& 3D & 0.0462 & 0.107 & 0.113 & -- & -- \\
\midrule
\multirow{2}{*}{Stress ($\gamma =1 $)} & 2D & 0.118 & 0.131 & 0.0814 & 0.123 & 0.0745 \\
& 3D & 0.0384 & 0.0928 & 0.0930 & -- & -- \\
\midrule
\end{tabular}
\caption{Tuned $\zeta$ values.}
\label{tab:zeta}
\end{table}

As mentioned, $\zeta \approx 0.1$ is a good initial choice for balancing the relative weights of the viscous, pressure gradient, and divergence operator sub-blocks of $A$. This value can be tuned so as to maximise multigrid performance. To this end, we considered a prototype model problem consisting of a constant-viscosity single-phase Stokes problem with periodic boundary conditions on a uniform $n \times n\,(\times\,n)$ grid; through numerical experiment, we then examined $\rho$ as a function of $\zeta$, where $\rho$ is the multigrid convergence rate using an \SAI1 smoother. For a fixed grid size $n$, we observed that $\rho$ is a smooth convex function of $\zeta$; a simple golden section search algorithm was implemented to locate the minimum and we executed this numerical search on progressively finer grids until the optimal value stabilises. This search was performed in 2D and 3D and on both forms of the Stokes equation as well as a variety of polynomial degrees $\polydeg$; \Cref{tab:zeta} contains the results of this study.\footnote{These values may depend on the specific choice of polynomial basis being used in the LDG framework; the values in \cref{tab:zeta} correspond to a tensor-product Legendre polynomial basis.} The values listed in \cref{tab:zeta} are used in the following results.

\subsection{Numerical experiments}
\label{sec:stokesexp}

In this section we examine multigrid performance for a variety of smoother and solver combinations, on a variety of Stokes problems, for a variety of polynomial degrees $\polydeg$, in two and three dimensions, and across multiple grid sizes. Our analysis uses the same methodology as that used for elliptic interface problems, see \cref{sec:scalarexp}. In particular, we apply the same V-cycle of \cref{algo:vcycle}: the interpolation and restriction operators are analogous to that for scalar elliptic problems, i.e., coarse mesh corrections are transferred to the fine mesh via polynomial injection, while fine mesh residuals are transferred to the coarse mesh via $L^2$ projection; equivalently, they can be defined via their scalar counterparts acting on the individual components of the velocity and pressure field. However, compared to the scalar case, we consider a reduced set of smoother/solver combinations. In particular, the Stokes operator matrix $A$ and corresponding V-cycle $V$ are both symmetric indefinite, so we cannot use CG as a backend solver. Consequently we need only explore combinations of (left-preconditioned) GMRES coupled with V-cycles that use block Jacobi, block Gauss-Seidel (with various combinations of pre- and post-smoother sweeping order), and block SAI (with various combinations of using \cref{eq:B1} or \cref{eq:BT} as pre- and post-smoother).

\subsubsection{Single-phase, steady-state, standard-form Stokes problems}
\label{sec:stokesexp1}

\begin{figure}%
\centering%
\includegraphics[scale=\orbscale]{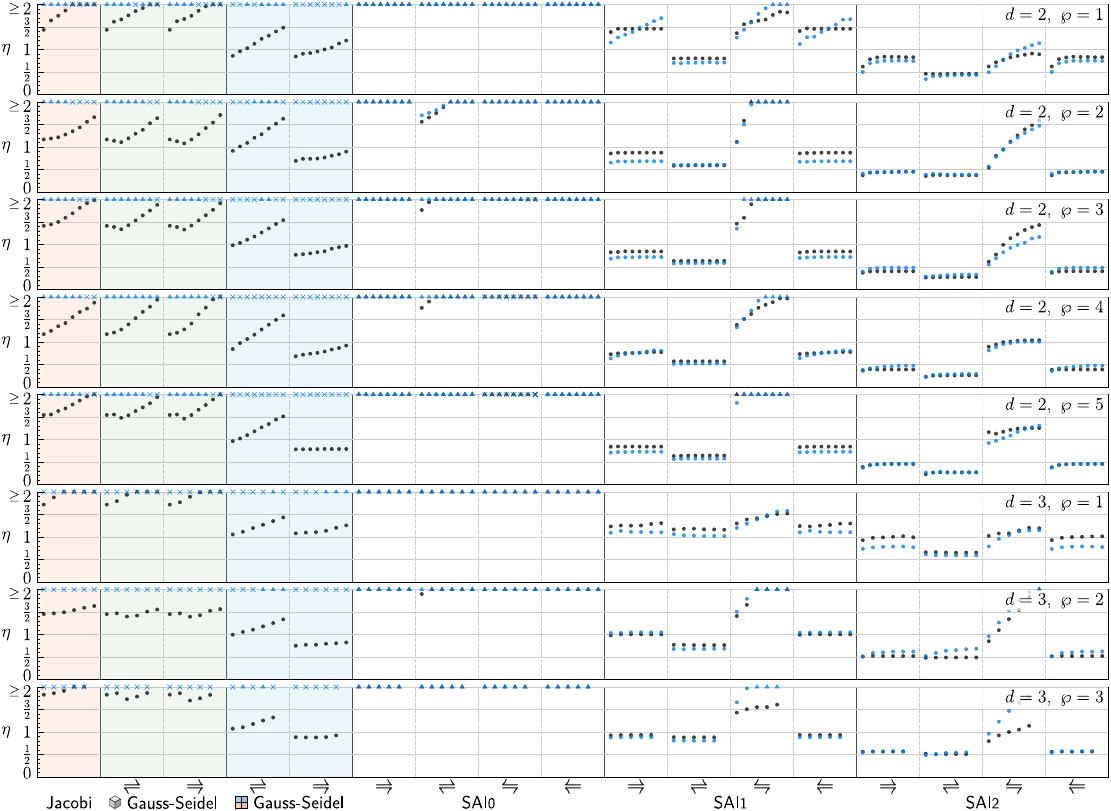}%
\caption{Multigrid solver performance for the single-phase steady-state standard-form Stokes problem considered in \cref{sec:stokesexp1}, with periodic boundary conditions. \seelegend{} Black (resp., blue) symbols correspond to the case when pressure penalty stabilisation is enabled (resp., disabled).}%
\label{fig:06}%
\end{figure}

\textbf{Periodic boundary conditions.} We begin with perhaps the simplest case, that of a single-phase steady-state standard-form Stokes problem with $\mu \equiv 1$ on the unit rectangular domain $\Omega = (0,1)^d$, applying periodic boundary conditions and uniform Cartesian grids. We consider two scenarios: (i) LDG pressure stabilisation is applied, using the optimised penalty parameters computed and explored in detail in \cite{flame}; and (ii) applying no LDG pressure stabilisation whatsoever, effectively so that the $E_p$ term in \cref{eq:blockform} is zero. We include the latter possibility mainly as a point of interest, as SAI-based smoothers can sometimes lead to faster solvers when there is no pressure stabilisation. \cref{fig:06} collects the numerical results of this study. Several observations can be made:
\begin{itemize}
    \item Out of the various block Gauss-Seidel smoother options, \xGSOOGL{} is markedly the fastest, provided we include LDG pressure stabilisation. On the other hand, processor-block Gauss-Seidel is as slow as damped Jacobi. As such, at least for Stokes problems using Gauss-Seidel solvers, it is worthwhile to use a fully undamped multi-coloured approach.
    \item In contrast to elliptic interface problems, \SAI{0}-based methods for Stokes are clearly too slow for practical use. This finding holds universally across all of the Stokes problems we tested; consequently, we shall remove \SAI{0} from consideration in the remainder of this article.
    \item The best performing \SAI1 method is that of \xSaiOIGL[1]{} and is faster than the best performing \xGS{} method. For example, with $d = \polydeg = 2$, \xSaiOIGL[1]{} yields $\eta \approx 0.6$ whereas \xGSOOGL{} approaches $\eta \approx 0.9$, representing a $1.5$-fold reduction in GMRES iteration count for the former. By this metric alone, this improves upon prior results \cite{flame}. 
    \item As for scalar elliptic interface problems, \xSaiIOGL[1]{} and \xSaiIOGL[2]{} stand out as being markedly ineffective. In other words, we see that it is suboptimal to use \cref{eq:BT} as a pre-smoother and \cref{eq:B1} as a post-smoother. The same is true for the remainder of the test problems in this paper.
    \item In regards to disabling LDG pressure stabilisation, this clearly renders every Jacobi- and Gauss-Seidel-based solver entirely ineffective. On the other hand, disabling pressure stabilisation can mildly speed up some of the SAI-based solvers. This result suggests that SAI smoothers can be more robust to a wider range of user-defined LDG penalty parameter values. Despite this speedup, we emphasise that pressure stabilisation is usually necessary to guarantee wellposedness of the equal-degree LDG Stokes discretisation; roughly speaking, it is possible to disable it in the current setting owing to the use of uniform Cartesian grids, simplified boundary conditions, and one-sided fluxes. For the remainder of this paper, we keep LDG pressure stabilisation enabled, using the parameter choices discussed in \cite{flame}.
\end{itemize}
In summary, putting aside the obviously ineffective solver combinations, there are a variety of multigrid solvers that can match the convergence rates of classical Poisson-style geometric multigrid methods, an outcome that can be elusive for saddle point problems like Stokes \cite{BenziGolubLiesen2005}.

\begin{figure}%
\centering%
\includegraphics[scale=\orbscale]{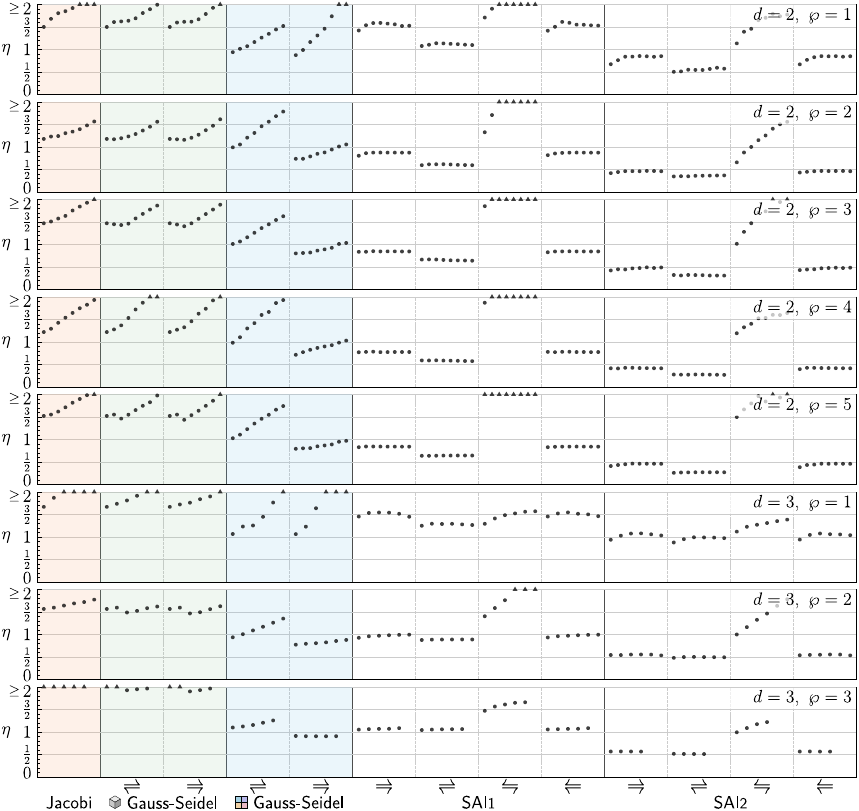}%
\caption{Multigrid solver performance for the single-phase steady-state standard-form Stokes problem considered in \cref{sec:stokesexp1}, with velocity Dirichlet boundary conditions. \seelegend{}}%
\label{fig:07}%
\end{figure}

\textbf{Dirichlet boundary conditions.} In addition to periodic boundary conditions, we also tested the case of Dirichlet boundary conditions on the velocity field (on the same single-phase, steady-state, standard-form Stokes problem). The corresponding numerical results are shown in \cref{fig:07}. Overall, they are similar to the case of periodic boundary conditions, except for the situation when $\polydeg = 1$, i.e., LDG discretisations using bilinear/trilinear polynomials. In this scenario we observe the best performing Gauss-Seidel method, \xGSOOGL{}, decays as the grid is refined, in both 2D and 3D. The same observation was made in the author's prior work \cite{flame} and currently lacks an explanation; at the same time we note that $\polydeg = 1$ is generally not of practical interest in the context of high-order LDG methods. Curiously, however, the best performing SAI methods on $\polydeg = 1$ are not subject to this slow degradation---their GMRES iteration counts remain constant as the mesh is refined.

\subsubsection{Single-phase, steady-state, stress-form Stokes problems} 
\label{sec:stokesexp2}

\begin{figure}%
\centering%
\includegraphics[scale=\orbscale]{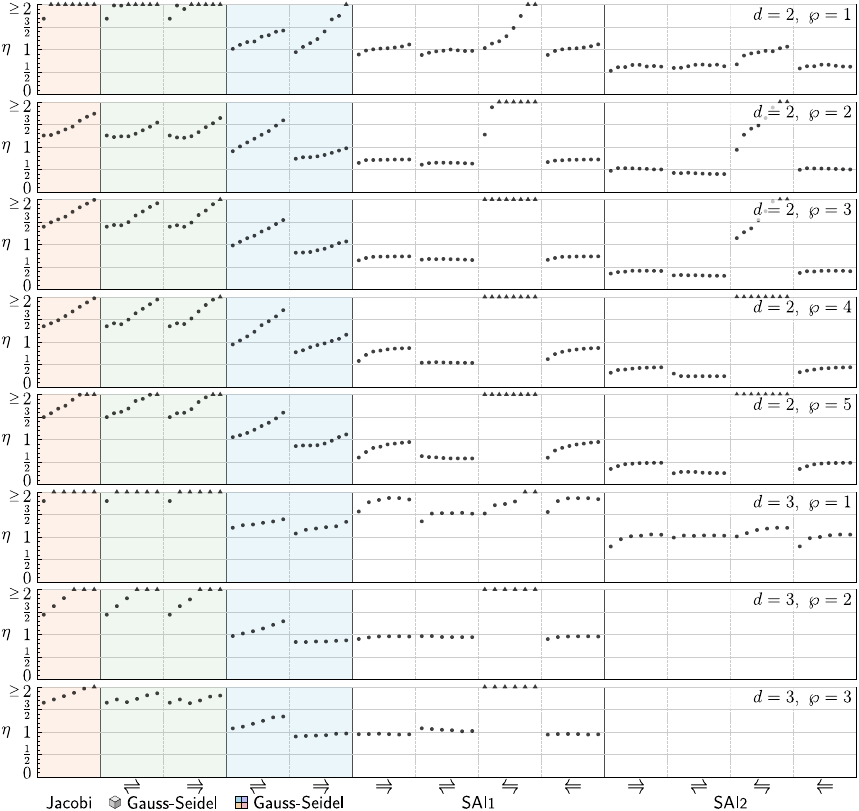}%
\caption{Multigrid solver performance for the single-phase steady-state stress-form Stokes problem considered in \cref{sec:stokesexp2}. \seelegend{}}%
\label{fig:08}%
\end{figure}

We next test performance of the various multigrid solvers when stress boundary conditions are imposed, i.e., ${\bm \sigma} \cdot \vn = \vh_\partial$ on $\partial \Omega$, in which case the appropriate form of the Stokes equations is the stress form ($\gamma = 1$). We consider again a single-phase problem on $\Omega = (0,1)^d$ with $\mu \equiv 1$ and uniform Cartesian grids. \cref{fig:08} compiles the numerical results;\footnote{On three-dimensional Stokes problems, \SAI2-based solvers are significantly more expensive to build, chiefly because of the increased stencil sizes in 3D and the cubic scaling in block size. Accordingly, owing to limited computing resources, some of the results shown in \cref{fig:08,fig:09,fig:10,fig:11} forgo testing of \SAI2-based multigrid solvers.} in large part they are similar to earlier results, showing that the various solvers perform similarly on both standard form and stress form Stokes problems. Once again, the best performing \SAI1 method is that of \xSaiOIGL[1]{} and is often faster than the best performing \xGS{} method. For example, with $d = \polydeg = 2$, \xSaiOIGL[1]{} yields $\eta \approx 0.6$ whereas \xGSOOGL{} yields $\eta$ approaching $\approx 1$, representing a $1.6$-fold reduction in GMRES iteration count for the former.

\subsubsection{Multiphase, steady-state, stress-form Stokes problems}
\label{sec:stokesexp3}

\begin{figure}%
\centering%
\includegraphics[scale=\orbscale]{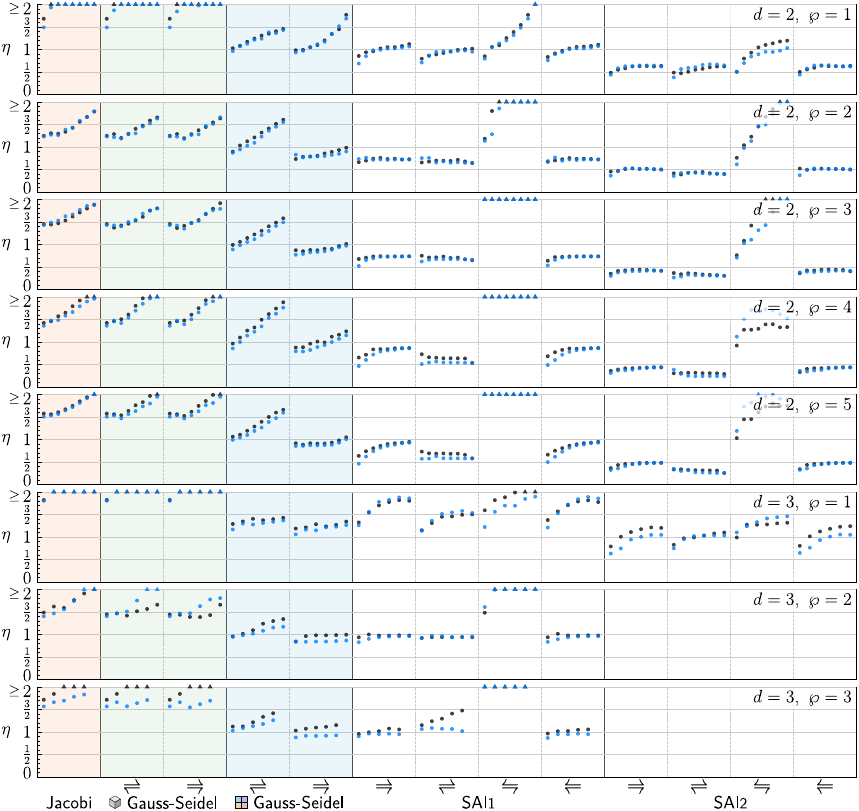}%
\caption{Multigrid solver performance for the multiphase steady-state stress-form Stokes  problem considered in \cref{sec:stokesexp3}. \seelegend{} Black (resp., blue) symbols correspond to the viscosity ratio $\mu_0/\mu_1 = 10^{-4}$ (resp., $10^4$). Other viscosity ratios, including $10^{\pm2}$ and $10^{\pm6}$, yield nearly identical results.}%
\label{fig:09}%
\end{figure}

In this test we consider a substantially more challenging case of a two-phase Stokes problem wherein the viscosity exhibits a jump multiple orders in magnitude across an embedded interface. Let $\Omega = (0,1)^d$ be divided into an interior rectangular phase $\Omega_0 = (\tfrac14, \tfrac34)^d$ and an exterior phase $\Omega_1 = \Omega \setminus \overline{\Omega_0}$, each with viscosity $\mu_i$. Two different viscosity ratios are evaluated: $\mu_0/\mu_1 = 10^4$ and $\mu_0/\mu_1 = 10^{-4}$. We consider the stress form of the multiphase Stokes equations \cref{eq:govern1,eq:govern2,eq:govern3} for which interfacial jump conditions are imposed on both velocity and stress on $\Gamma = \partial \Omega_0 \cap \partial \Omega_0$, and, for simplicity, periodic boundary conditions on $\partial \Omega$. \cref{fig:09} compiles the multigrid speed results for the two viscosity ratios. As in other results, for each $\polydeg$ in each dimension $d$, there is at least one smoother choice that results in an effective multigrid solver. This favourable outcome is attributed in part to the use of viscosity-upwinded numerical fluxes, which bias the LDG discretisation to the mutual benefits of multigrid solver speed and high-order accuracy \cite{flame}, and in the case of SAI-based methods, the use of balancing as discussed in \cref{sec:stokesbalance}. Besides the viscosity ratios of $10^{\pm4}$, we also tested other ratios including $10^{\pm2}$ and $10^{\pm6}$; essentially identical results are seen in these cases.

\subsubsection{Unsteady Stokes problems}
\label{sec:stokesexp4}

So far we have focused on steady-state Stokes problems; our last set of test cases examine unsteady/time-dependent/nonstationary Stokes problems. The governing equations consist of keeping the interfacial jump conditions \cref{eq:govern2} and boundary conditions \cref{eq:govern3} but amending \cref{eq:govern1} to include a density term that would arise from, e.g, a time-stepping method applied to the incompressible Navier-Stokes equations whereby the advection term is treated explicitly and the viscous term implicitly. The amended form reads
\begin{equation} \label{eq:govern1unsteady} \left. \begin{aligned}  \frac{\rho_i}{\delta} \vu -\nabla \cdot \bigl(\mu_i (\nabla \vu + \gamma \nabla \vu^\trans) \bigr) + \nabla p &= \mathbf f \\ -\nabla \cdot \mathbf \vu &= f \end{aligned} \right\} \text{ in } \Omega_i. \end{equation}
Here, $\rho_i > 0$ is a phase-dependent density while $\delta > 0$ is a parameter proportional to the time step $\Delta t$ of a temporal integration method. The corresponding LDG discretisation is a straightforward modification to the steady-state problem and leads to
\begin{equation} \label{eq:blockformunsteady} \begin{pmatrix} \tfrac{1}{\delta} M_\rho + \mathscr{A} & \mathscr{G} \\ -\mathscr{D} & -E_p \end{pmatrix} \begin{pmatrix} \vu_h \\ p_h \end{pmatrix} = \begin{pmatrix} {\mathbf b}_{\vu} \\ b_p \end{pmatrix}\!, \end{equation}
where $M_\rho$ is a $\rho$-weighted mass matrix. Beyond this simple modification, what makes the unsteady problem subtle is the system's dependence on the relative strengths of the viscous operator and the newly-added density-weighted term. When $\mu \delta/\rho$ is sufficiently large, the dominant operator is the steady-state Stokes equations, with a small $\rho/\delta$-weighted identity shift added to the viscous operator; in this case, a good steady-state Stokes solver can be effective. Conversely, when the viscous effects are weak and $\rho/(\delta \mu)$ is sufficiently large, e.g., when the mesh is unable to resolve viscous effects, then \cref{eq:blockformunsteady} approximately reduces to solving a kind of Helmholtz decomposition, closely connected to Chorin's projection method for incompressible fluid dynamics \cite{Chorin}. The role of the pressure variable $p$ changes in these two extremes, from being a Lagrange multiplier enforcing a divergence constraint on $\vu$ for the former, to being the primary solution variable in an elliptic interface problem for the latter. Consequently, we may expect altered multigrid smoothing characteristics. For example, on a relatively deep multigrid hierarchy, it may be the case that viscous effects dominate on the finest grid, but become relatively weak on the coarsest grid, so the balance can change across the hierarchy. In the author's prior work, see \cite{flame}, a simple method was developed to define the LDG pressure penalty parameter so as to effectively blend between these two extremes. In combination with a multi-coloured Gauss-Seidel method, the approach was found therein to be effective; we use the same approach here and compare the Gauss-Seidel approach with SAI-based smoothers.

\begin{figure}%
\centering%
\includegraphics[scale=\orbscale]{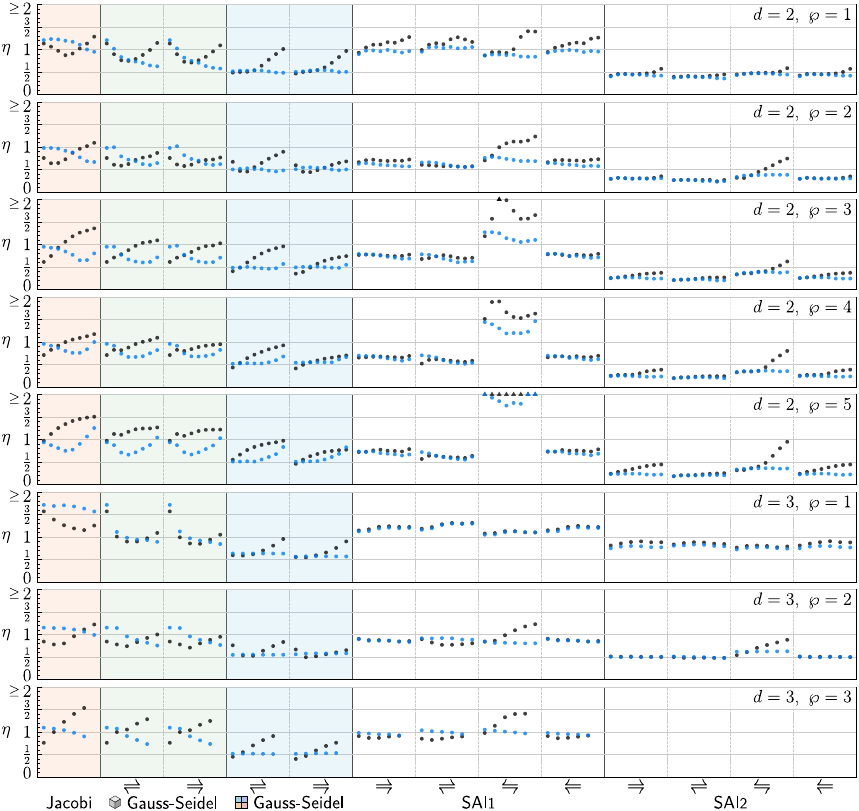}%
\caption{Multigrid solver performance for the single-phase unsteady standard-form Stokes problem considered in \cref{sec:stokesexp4}. \seelegend{} Black (resp., blue) symbols correspond to an effective Reynolds number of $10^2$ (resp., $10^4$).}%
\label{fig:10}%
\end{figure}

\textbf{Single-phase, unsteady, standard-form Stokes problems.} We consider two effective Reynolds numbers $\textsf{Re} = \rho U L/\mu$: $\textsf{Re} = 10^2$, representing a viscous-dominated case (but where the time-derivative operator nevertheless influences performance characteristics) and $\textsf{Re} = 10^4$, representing a case wherein the time-derivative operator definitively dominates. In both scenarios, the velocity and length scales are unitary, $U = 1$ and $L = 1$, while density is set to $\rho = 1$. In addition, we define $\delta = 0.1 h$, where $h$ is the element size of the primary-level uniform Cartesian grid, representing a typical scenario of applying the unsteady Stokes equations in a temporal integration method with unit-order CFL. For simplicity, periodic boundary conditions are applied. With this setup, \cref{fig:10} compiles the numerical results, measuring multigrid solver speed across a variety of solver choices. For \xGS{}-based solvers, in the case when $\mu = 10^{-2}$ we sometimes see a mild upwards trend in $\eta$ as the primary mesh is refined: the same behaviour was observed in prior work \cite{flame} and is attributed to the fact that as the grid is refined, eventually we obtain a slightly perturbed unsteady problem and thus recover the multigrid solver speeds seen in \cref{fig:06}. On the other hand, SAI-based methods (excluding the outliers \xSaiIOGL[1]{} and \xSaiIOGL[2]{}) appear to be less dependent on both mesh resolution and effective Reynolds number, having relatively constant convergence rates across all grid sizes.

\begin{figure}%
\centering%
\includegraphics[scale=\orbscale]{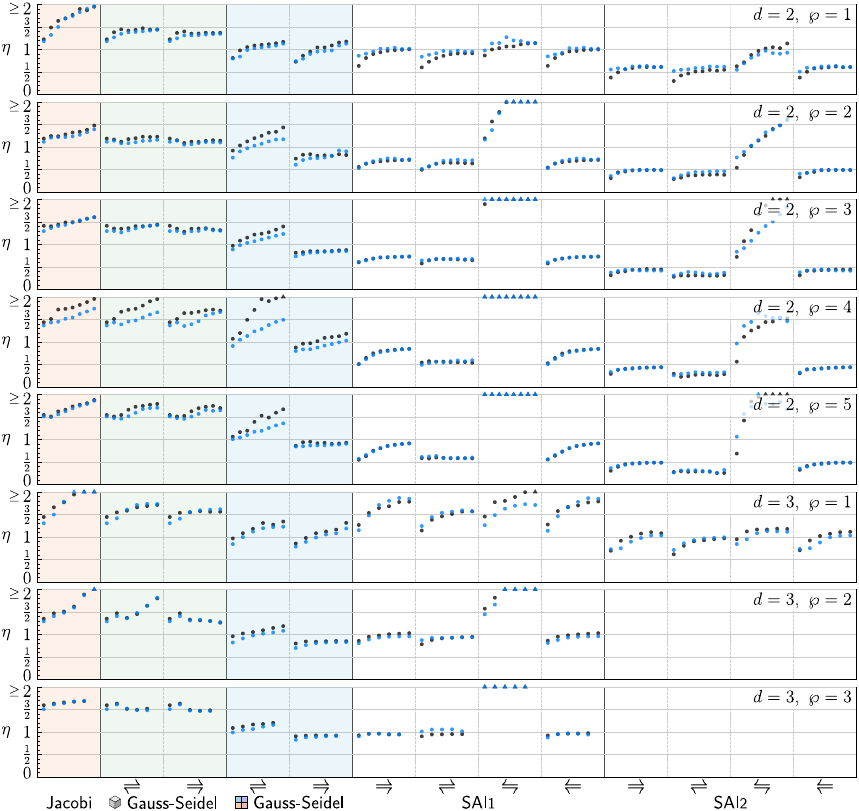}%
\caption{Multigrid solver performance for the multiphase unsteady stress-form Stokes problem considered in \cref{sec:stokesexp4}. \seelegend{} Black (resp., blue) symbols correspond to a gas droplet surrounded by water (resp., water droplet surrounded by gas).}%
\label{fig:11}%
\end{figure}

\textbf{Multiphase, unsteady, stress-form Stokes problems.} Finally, in our last two examples we combine the challenging aspects of a two-phase problem in which viscosity and density both jump by several orders in magnitude across an embedded interface, together with the effects of the newly-added time-dependent term. We consider two scenarios: a water bubble surrounded by gas, and a gas bubble surrounded by water. Specifically, $\rho_{\text{water}} = 1$, $\rho_{\text{gas}} = 0.001$, $\mu_{\text{water}} = 1$, and $\mu_{\text{gas}} = 0.0002$ (approximately accurate values for water and air at ambient temperature in CGS units). Here, the ``bubble'' refers to an interior rectangular phase $\Omega_0 = (\tfrac14,\tfrac34)^d$ while the exterior phase is given by $\Omega_1 = \Omega \setminus \overline{\Omega_0}$ where $\Omega = (0,1)^d$. We solve \cref{eq:govern1unsteady,eq:govern2,eq:govern3} in stress form, $\gamma = 1$, with stress jump conditions across the water-gas interface, and velocity Dirichlet boundary conditions. As before, we set $\delta = 0.1 h$ to represent, say, solving the incompressible Navier-Stokes equations through explicit treatment of the advection term and implicit treatment of the viscous term. \cref{fig:11} compiles the numerical results on this multiphase unsteady Stokes problem. Besides $\polydeg = 1$, which as previously noted is a kind of outlier, we observe the best {\xGS}-, {\SAI1}- and {\SAI2}{}-based solvers all exhibit the performance of classical multigrid methods. For example, to reduce the residual by a factor of $10^8$ in the case of $d=2$ and $\polydeg = 5$, each of the corresponding solvers need around 7.2, 4.8, and 2.4 GMRES iterations, respectively, essentially independent of the grid size.

\section{Concluding Remarks}
\label{sec:conclusions}

In this study, we explored a variety of multigrid methods for solving elliptic interface problems and multiphase Stokes problems discretised through high-order LDG methods. Our approach used a standard multigrid V-cycle and tested a variety of element-wise block smoothers: block Jacobi, block Gauss-Seidel, and block SAI methods. In addition, a simple balancing method was developed to effectively and automatically treat: (i) for multiphase problems, ellipticity/viscosity coefficients varying over multiple orders of magnitude; and/or (ii) for Stokes problems, the viscosity- and mesh size-dependent scaling of the viscous, pressure gradient, and velocity divergence operators. Across an array of 2D and 3D test cases, with different choices of polynomial degrees $\polydeg$ and on multiple grid sizes, we examined a variety of multigrid smoother and solver combinations and found that various choices yield convergence rates matching that of classical geometric multigrid methods; e.g., in many cases, a $10^{10}$-fold reduction in the residual requires 4 to 8 iterations.

In the case of elliptic interface problems, the balancing method can be viewed as a simple diagonal preconditioning, or, alternatively, as a modification to the standard SAI construction whereby the Frobenius norm is replaced by a weighted Frobenius norm. The same is true for Stokes problems, but in that case the balancing also serves another role, which is to balance the interplay of the velocity and pressure degrees of freedom. Within each diagonal block of the discrete Stokes operator, the balancing method rescales the viscous operator, pressure gradient, and velocity divergence block components to have Frobenius norm $1$, $\zeta$, and $\zeta$, respectively. Here, $\zeta$ is a user-defined parameter that controls the weighting in the SAI's approximate inverse. Intuitively, one can choose $\zeta$ so that the high frequency error components of velocity and pressure are smoothed at the same rate. We found that $\zeta \approx 0.1$ works well across all tested Stokes problems, in 2D and 3D, and across polynomial degrees $\polydeg$, though $\zeta$ can be tuned further, per \cref{tab:zeta}.

\begin{table}%
\centering%
\sffamily\footnotesize%
\newcommand{\ssize}[2]{$#1 \times #2$}%
\begin{tabular}{lcc|ccc|ccc}
& & \multicolumn{1}{c}{Block} & \multicolumn{3}{c}{SAI stencil size in blocks} & \multicolumn{3}{c}{SAI smoother  cost}\\
Problem & $d$ & dimension & \SAI0 & \SAI1 & \SAI2 &\SAI0 & \SAI1 & \SAI2 \\
\midrule
\multirow{2}{*}{Poisson}
& 2D & $(\polydeg + 1)^2$ & \ssize{5}{1} & \ssize{13}{5} & \ssize{25}{13} & 6 & 10 & 18 \\ 
& 3D & $(\polydeg + 1)^3$ & \ssize{7}{1} & \ssize{25}{7} & \ssize{63}{25} & 8 & 14 & 32 \\
\midrule
\multirow{2}{*}{Stokes, standard form}
& 2D & $3(\polydeg + 1)^2$ & \ssize{5}{1} & \ssize{13}{5} & \ssize{25}{13} & 6 & 10 & 18 \\ 
& 3D & $4(\polydeg + 1)^3$ & \ssize{7}{1} & \ssize{25}{7} & \ssize{63}{25} & 8 & 14 & 32 \\
\midrule
\multirow{2}{*}{Stokes, stress form}
& 2D & $3(\polydeg + 1)^2$ & \ssize{7}{1} & \ssize{19}{7} & \ssize{37}{19} & 8 & 14 & 26 \\ 
& 3D & $4(\polydeg + 1)^3$ & \ssize{13}{1} & \ssize{55}{13} & \ssize{147}{55} & 14 & 26 & 68 \\
\midrule
\end{tabular}
\caption{Individual block dimensions, SAI stencil sizes (in blocks), and SAI smoother cost (per element, in block units), for canonical Poisson and Stokes problems on uniform $d$-dimensional Cartesian grids using tensor-product polynomials of one-dimensional degree $\polydeg$.}
\label{tab:sizes}
\end{table}

We have yet to discuss the associated computational cost of SAI-based methods. This can be divided into two parts: (i) building the SAI itself, and (ii) once (pre)computed, applying the SAI as a smoother via \cref{eq:B1} or \cref{eq:BT}.
\begin{enumerate}
    \item Each block column of $B_\ell$ can be computed independently of (and possibly concurrently with) every other block column; each calculation leads to an overdetermined least squares problem whose corresponding matrix has row count equal to $\mathsf{m}$ blocks and column count equal to $\mathsf{n}$ blocks, where $\mathsf{m} \times \mathsf{n}$ is the stencil size. \Cref{tab:sizes} lists the stencil sizes for the canonical\footnote{Here, ``canonical'' refers to the basic case arising from a uniform Cartesian grid, ignoring the influence of interfacial jump and boundary conditions; e.g., the canonical Poisson problem leads to a 5-point or 7-point Laplacian stencil in 2D and 3D, respectively.} Poisson, standard-form Stokes, and stress-form Stokes problems. Using a standard dense least squares method (e.g., QR) the corresponding solver cost is $\mathcal{O}\bigl( (\polydeg + 1)^{3d} \mathsf{m} \mathsf{n}^2 \bigr)$, where the prefactor accounts for the cubic scaling in the block size. Obviously, this can be substantial, especially for the large stencil sizes encountered in 3D.

    Fortunately, it is often the case that many block columns of $B_\ell$ have the \textit{same} block-sparse solution, up to a one-to-one mapping of the elements involved. Therefore, a representative solution can be computed once, cached, and reused across all applicable columns, and even across the multigrid hierarchy itself. Indeed, the blockwise stencil of the (balanced) primary discrete Laplacian/Stokes operator $A$ is often the same on most levels (at least, on well-structured grids, away from interfaces and boundaries); therefore, columns of $B_\ell$ can be reused across levels as well. This can lead to a substantial speedup; e.g., on a uniform Cartesian grid with periodic boundary conditions, just a handful of SAI columns need computation, one invoked by the fine mesh, and a few corresponding to the bottom-most levels where the stencil of $A$ might wrap around and overlap. We used this approach in this study: our cacheing mechanism stored the computed SAI block-column solutions in a hash table combined with one-to-one mapping algorithm to match the stencil structure and blocks involved.

    On the other hand, the potential speedups of cacheing diminish whenever there is insufficient repetition within $A$, e.g., as would arise when the domain boundary or interfaces are sufficiently dense, or the mesh sufficiently unstructured mesh or is too frequently changing structure in a dynamical setting, or when the blocks or stencil sizes themselves are too large, etc. As an extreme example, for a stress-form Stokes problem, a single block column calculation for \SAI2 typically involves a $147 \times 55$ block stencil; for $\polydeg = 3$, each block is of size $256 \times 256$, leading to a least squares system of size about $37632$ by $14080$ with corresponding solver flop count around $10^{13}$; clearly prohibitive to solve any more than a handful of times (on a single core). For many applications, e.g., on unstructured or dynamic meshes, or anytime where the precomputation of SAI cannot be amortised, \SAI2 methods are likely too costly, and potentially \SAI1 as well.

    \item Regarding the cost of applying the SAI smoother, one iteration corresponds to either the update $\xout = \xin - B(A \xin - b)$ or the update $\xout = \xin - B^\trans(A \xin - b)$. In both cases, this amounts to a block-sparse residual calculation followed by a block-sparse matrix-vector multiplication. For the canonical problems mentioned above, \cref{tab:sizes} indicates the predominant cost per mesh element in block units. (For example, on a 2D Poisson problem, $A$ has a stencil size of 5, so \SAI0 (resp., \SAI1, \SAI2) costs $5 + 1$ (resp., $5 + 5$, $5 + 13$) block matrix-vector operations per element.) In particular, note that \SAI1 has the same stencil size as $A$, therefore one smoother application of \SAI1 is roughly twice the cost of one smoother application of either block Jacobi or block Gauss-Seidel. Consequently, \SAI1 is perhaps beneficial only when its convergence rate is appreciably faster than Gauss-Seidel's, say, so that it leads to sufficiently fewer V-cycle iterations or the ability to decrease the number of pre- and post-smoother steps. Likewise, if both \SAI1 and \SAI2 result in good multigrid speeds, then \SAI2 would only be competitive if its convergence rate is significantly better than \SAI1: referring to \cref{tab:sizes}, the per-element cost of applying \SAI2 is about twice that of \SAI1; on the other hand, for most of the results presented here, the best \SAI2 solver has $\eta$ at most two-times (sometimes three-times) smaller than that of the best \SAI1 solver, so in many cases \SAI2 is not worthwhile once we factor in its construction overhead. Naturally, the trade-off between Gauss-Seidel vs.~\SAI1 vs.~\SAI2, improvements in CG/GMRES iteration counts, changing the number of pre- and post-smoothing steps, etc., all depend on a variety of implementation aspects, such as the use of high-arithmetic density operations like level-3 BLAS, vectorisation, memory locality and distributed vs.~shared memory, V-cycle communication/synchronisation costs, and so forth.
\end{enumerate}

Beyond the findings already discussed in the numerical experiments of \cref{sec:scalarexp} and \cref{sec:stokesexp}, some additional conclusions can be made:
\begin{itemize}
    \item In every single test problem, there was a solver based on \SAI1 and \SAI2 that led to good multigrid behaviour, i.e., V-cycle iteration counts that remain small and bounded as the primary grid is refined. \SAI2 yields fewer iterations than \SAI1, but not enough to offset its comparatively higher construction costs, as discussed above; in other words, of the two approaches, \SAI1 is the more computationally efficient solver in terms of minimising solver time. Moreover, in every single test problem, \SAI0 is always slower than damped Jacobi, and in many cases (notably Stokes problems) impractically slow; consequently element-blockwise \SAI0-based smoothers appear to have no compelling utility, at least for LDG problems.
    \item Across a significant majority of test problems, multi-coloured Gauss-Seidel yielded highly efficient multigrid solvers. Oftentimes its corresponding $\eta$ is similar to the best \SAI1-based solver, however the latter has comparatively higher computational cost. Therefore, in many settings, multi-coloured Gauss-Seidel leads to the fastest-possible solver measured by computation time. However, there are some instances where multi-coloured Gauss-Seidel is suboptimal because $\eta$ slowly grows as the grid is refined (though it will eventually stabilise). These cases correspond to bilinear/trilinear polynomials in Stokes problems, i.e., when $\polydeg = 1$ in \cref{fig:06,fig:07,fig:08,fig:09}, and perhaps are not of interest in the context of high-order LDG methods. Nevertheless, \SAI1 is competitive for $\polydeg = 1$ as it results in small, relatively constant iteration counts.
    \item In addition, multi-coloured Gauss-Seidel may sometimes be too complex to implement because of the need for a good colouring, which for some applications can be intricate to determine, e.g., on high resolution 3D unstructured meshes. As such, it is often simpler to implement a processor-block version of Gauss-Seidel \cite{AdamsBrezinaHuTuminaro}, but this usually requires damping, that may also be intricate to tune. In a majority of test problems, the fastest \SAI1-based solver can be substantially faster than processor-block Gauss-Seidel; moreover, SAI is more parallel friendly (in fact, essentially ``embarrassingly parallel'') and does not require damping, so in these instances, \SAI1 may again be competitive.
    \item We also note that SAI-based smoothers can be less sensitive to the particular choice of LDG penalty stabilisation parameters. In prior work \cite{flame} it was found that for Stokes problems, Gauss-Seidel works best when the pressure penalty value is somewhat finely-tuned. On the other hand, some of our results (see \cref{fig:06}) show that SAI can work even in the degenerate case of no pressure stabilisation whatsoever. Meanwhile, the penalty parameters associated with weakly enforcing Dirichlet boundary conditions and interfacial jump conditions can also affect the convergence rates of V-cycles \cite{dgmg}. Using a V-cycle as a preconditioner in a Krylov method can mitigate these effects, but the issue of choosing good penalty parameters remain. Some of our numerical experiments (not shown here) indicated that SAI-based smoothers are robust to a broader range of penalty parameter choices.
    \item In the numerical results we saw that the pre- and post-smoother ordering plays a crucial role in attaining fast multigrid convergence. The fastest combinations often result in an asymmetric V-cycle and thus cannot be used with CG; however, the speedup gained more than offsets the increased memory footprint of GMRES. If CG is abandoned as the primary Krylov solver and symmetry requirements on the V-cycle are dropped, then various additional possibilities open up for further multigrid optimisation, including V-cycles or W-cycles with differing pre- and post-smoother counts, combining different smoothers, and nonsymmetric balancing/preconditioning approaches that solely pre-scale or post-scale the system.
    \item We examined here a high-order LDG framework using a tensor-product polynomial basis on Cartesian grids, quadtrees, and octrees. Nevertheless, it is expected many of our findings would carry over to other basis choices as well as more unstructured meshes and to other kinds of DG methods. For example, we also tested SAI-based multigrid methods on LDG methods that use total-degree polynomial bases\footnote{For example, $\polydeg = 2$ in 2D would correspond to $\spanop \{1,x,y,x^2,xy,y^2\}$.}; similar multigrid behaviour was observed. Preliminary work on combining SAI methods with LDG on implicitly defined meshes \cite{ImplicitMeshPartOne,ImplicitMeshPartTwo} also shows promise, particularly in regard to robustly and efficiently treating complex geometry and arbitrarily small cut cells.
\end{itemize}

Finally, we mention some other possibilities for future research. We considered here element-wise block smoothers and SAI block-sparsity patterns based on the powers of $A$; perhaps other strategies could be useful. For example, in a nodal DG method, one could shrink the sparsity pattern of \SAI1 by including only the central element and the shared nodes of its neighbours. For Stokes problems, another possibility may be to block the velocity and pressure variables separately and similarly refine the sparsity patterns of SAI. Meanwhile, SAI methods have been shown to be effective for anisotropic finite difference and finite volumes problems \cite{1997doi:10.1137/S1064827594276552,2000doi:10.1137/S106482759833913X,2000doi:10.1137/S0895479899339342,2001doi:10.1137/S1064827500380623,2002BROKER200261,2010doi.org/10.1155/2010/930218}; the same likely holds for high-order LDG methods applied to anisotropic problems. Last, all of these findings could also be useful to inform the design of algebraic multigrid methods \cite{2002BROKER200261}.

\section*{Acknowledgements}

This work was supported in part by the U.S. Department of Energy, Office of Science, Office of Advanced Scientific Computing Research's Applied Mathematics Competitive Portfolios program under Contract No.~DE-AC02-05CH11231, and by a U.S. Department of Energy, Office of Science, Early Career Research Program award. Some computations used resources of the National Energy Research Scientific Computing Center (NERSC), a U.S. Department of Energy, Office of Science User Facility operated at Lawrence Berkeley National Laboratory.

\bibliographystyle{siamplain}
\bibliography{references}

\end{document}